\documentclass{amsart}

\usepackage{amssymb,amscd,xypic,amsmath,graphics}

\newtheorem{theorem}{Theorem}[section]
\newtheorem{lemma}[theorem]{Lemma}
\newtheorem{proposition}[theorem]{Proposition}
\newtheorem{corollary}[theorem]{Corollary}

\theoremstyle{definition}

\newtheorem{definition}[theorem]{Definition}
\newtheorem{remark}[theorem]{Remark}

\newtheorem{example}[theorem]{Example}

\newcommand\bk{{\Bbbk}}

\newcommand\bP{{\mathbb P}}
\newcommand\bQ{{\mathbb Q}}
\newcommand\bN{{\mathbb N}}
\newcommand\bZ{{\mathbb Z}}
\newcommand\bR{{\mathbb R}}

\newcommand\bT{{\mathbb T}}

\newcommand\cO{{\mathcal O}}

\newcommand\cX{{\mathcal X}}
\newcommand\cB{{\mathcal B}}

\newcommand\nd{\mathop{\rm nd}\nolimits}
\newcommand\cc{\operatorname{cc}}
\newcommand\codim{\mathop{\rm codim}\nolimits}
\newcommand\Id{\mathop{\rm id}\nolimits}
\newcommand\GL{\mathop{\rm GL}\nolimits}
\newcommand\Hom{\mathop{\rm Hom}\nolimits}
\newcommand\opJ{\mathop{\rm J}\nolimits}
\newcommand\Aut{\mathop{\rm Aut}\nolimits}
\newcommand\supp{\mathop{\rm supp}\nolimits}
\newcommand\Ker{\mathop{\rm Ker}\nolimits}

\title{Self-dual projective toric varieties} 

\author[M. Bourel, A. Dickenstein and A. Rittatore]{Mathias Bourel, Alicia Dickenstein and Alvaro Rittatore}

\thanks{MB was partially supported by grant ``Variedades t\'oricas
proyectivas y dualidad'', CSIC, Uruguay;  AD was 
partially supported by UBACYT X064, CONICET PIP 112-200801-00483 and ANPCyT PICT 2008-0902, Argentina, 
and FCE 10018, Uruguay; AR was partially supported
by FCE 10018, Uruguay.} 

\begin{document}
\maketitle

\begin{abstract}

Let $T$ be  a torus over an algebraically closed field $\bk$
of characteristic $0$, and consider a projective  $T$-module $\mathbb P(V)$. 
We determine when a projective toric subvariety  $X\subset \mathbb
P(V)$ is self-dual, in terms of the configuration of weights of $V$.  
\end{abstract}

\section{Introduction}
\label{sec:introduction}

The notion of duality of projective varieties, which appears in various
branches of mathematics, has been a subject of study since the
beginnings of algebraic geometry \cite{kn:Gel,kn:tev05}.  Given an embedded
projective variety $X\subset \bP (V)$, its dual variety $X^*$ is the
closure in the dual projective space $\bP(V^\vee) $ of the hyperplanes 
intersecting  the regular points of $X$ non transversally.
 
A projective variety $X$ is {\em self-dual} if it is
isomorphic to its dual $X^*$ as embedded projective varieties.   The expected
codimension of the dual variety is  one. If
this is not the case, $X$ is said to be defective. 
Self-dual varieties other than hypersurfaces are defective varieties
with ``maximal'' defect.

Let  $\Bbbk$ be an algebraically closed field
of characteristic $0$. Let $T$ be an algebraic torus over $\Bbbk$  and
$V$ a finite dimensional $T$-module. In this paper 
we characterize self-dual projective toric varieties $X\subset
\bP(V)$ equivariantly embedded, in terms of the
combinatorics of the associated configuration of weights $A$ (cf. Theorems~\ref{teo:Bside}
and \ref{teo:Aside}) and in terms of the
interaction of the space of relations of these weights with the torus orbits (cf. 
Theorems~\ref{teoeqauto} and \ref{teo:constaut+gde}).
In particular, we show that $X$ is self-dual if and only if  $\dim X =
\dim X^*$ and  the smallest linear subspaces containing
$X=X_A$ and $X^*$ have the same dimension, see
Theorems~\ref{thm:autdim} and  \ref{thm:autdimgral}.

Given a basis of eigenvectors of $V$ and the configuration of weights
of the torus action on $V$, 
it is not difficult to check the equality of the dimensions of $X$ and
its dual (for instance, by means of the combinatorial characterization
of the tropicalization given in \cite{kn:DFS05}). 
But the complete classification of defective projective toric varieties
in an equivariant embedding is open in full generality and involves a complicated combinatorial
problem. For smooth toric varieties
this characterization is obtained in \cite{kn:DR03}; the case of $\bQ$-factorial toric
varieties is studied in \cite{kn:DRC05}.  For non necessarily normal
projective toric varieties of codimension two, a characterization is given in
\cite{kn:DS02}. This has been extended for codimensions three and four in \cite{kn:CC05}.

For smooth projective varieties,
a full list of self-dual varieties is known \cite{kn:Ein1,kn:Ein2,kn:tev05}. This list is indeed short
and reduces in the case of toric varieties  to hypersurfaces or Segre embeddings of $\bP^1 \times \bP^{m-1}$, for any
$m \geq 2$, under  the assumption that $\dim
X \leq \frac{ 2 \dim \bP(V)} {3}$. This 
was expected to be the whole classification under the validity of Hartshorne's 
conjecture \cite{kn:Ein1}.  We prove that this is indeed the whole list of self-dual 
smooth projective toric varieties in Theorem~\ref{th:regcase}. 

There exist some classical examples of self-dual non smooth varieties, as
the quartic Kummer surface. Popov and Tevelev gave new families of
non smooth self-dual varieties that come from actions of isotropy
groups of complex symmetric spaces on the projectivized nilpotent 
varieties of isotropy modules (\cite{kn:Pop02}, \cite{kn:PT04}). 
As a consequence of Theorem~\ref{teo:Bside}, it is easy to construct new families of
self-dual projective toric varieties in terms of the Gale dual 
configuration (see Definition~\ref{def:Gale}).   

A big class of self-dual toric varieties are the toric varieties associated to
Lawrence configurations (see Definition~\ref{df:Law}), which contain the configurations
associated to the Segre embeddings. Lawrence constructions are well known in the
domain of geometric combinatorics, where they
are one of the prominent tools  to visualize the geometry of
 	higher dimensional polytopes (see \cite[Chapter 6]{kn:Zi}); the commutative algebraic properties
 	of the associated toric ideals are studied in \cite{kn:bps}.
We show in Section~\ref{sec:examples}
other non Lawrence concrete examples for any
dimension bigger than $2$ and any codimension bigger than $1$. 

We also introduce the notion of strongly self-dual toric varieties
(see Definition~\ref{def:strong}),
which is not only related to the geometry of the
configuration of weights but also to number theoretic aspects. This concept
is useful for the study of the existence of rational multivariate hypergeometric
functions \cite{kn:Gel2,kn:cds01}.

In Section~\ref{sec:prelim} we gather some preliminary results about embedded
projective toric varieties and duality of projective varieties. In
Section~\ref{sec:orbits} we characterize self-dual  projective toric  varieties 
in terms of the geometry of the action of
the torus and we give precise assumptions under which self-dual projective
varieties are precisely those with maximal defect. 
In Section~\ref{sec:comb} we give two (equivalent) combinatorial
characterizations of self-duality. In Section~\ref{sec:examples} we collect several new
examples of self-dual (non smooth) projective toric varieties.
Finally, in Section~\ref{sec:strong} we study strongly self-dual toric varieties.

\bigskip
\noindent {\em Acknowledgments:}\/
We are grateful to Eduardo Cattani for his suggestion of the statement of 
Theorem~\ref{teo:Bside} and to Vic Reiner for his help in the statement
of Theorem~\ref{teo:Aside}. We also thank the  referee for posing helpful 
questions. 

\section{Preliminaries} \label{sec:prelim}

In this section we collect some well known results and useful
observations on projective toric 
varieties and duality of projective varieties.

\subsection{Actions of tori}

Let $T$ be an algebraic torus  over an algebraically closed field $\bk$
of characteristic $0$.  
We denote by $\cX(T)$ the lattice of characters  of $T$; recall that
$\Bbbk[T]=\bigoplus_{\lambda\in\cX(T)} \Bbbk \lambda$. 
Any finite dimensional
rational $T$-module $V$, $\dim V=n$, decomposes as
a direct sum of irreducible representations 
\begin{equation} \label{eq:decomp}
V\cong \bigoplus_{i=1}^n \Bbbk v_i,
\end{equation}
with $t\cdot v_i=\lambda_i(t)v_i$,   $\lambda_i\in \cX(T)$,  for all $t\in
T$.  

The action of $T$ on $V$ canonically induces an action $T\times
\bP(V)\to \bP(V)$ on the associated projective space, 
given by $t\cdot [v]=[t\cdot v]$, where $[v]\in \bP(V)$ denotes the
class of $v\in V \setminus \{0\}$. 
Recall that an irreducible $T$-variety $X$ is called {\em toric}\/ if there exists
$x_0\in X$ such that the orbit $\mathcal O(x_0)$ is open in $X$.

Let $A=\{\lambda_1,\dots,\lambda_n\}$ (which may contain repeated elements)
be the associated set of weights of a finite dimensional $T$-module
$V$ --- we call $A$ the {\em configuration of weights}\/ associated to the
$T$-module $V$. 
 To any basis  $\cB=\{v_1,\dots,v_n\}\subset
V$  of eigenvectors we can
associate a projective toric variety by
\[X_{V,\cB}=\overline{\mathcal
O\bigl(\bigl[{\textstyle\sum v_i}\bigr]\bigr)}\subset \bP(V).
\]

Denote  by $\bT^{n-1}=\bigl\{
\sum p_iv_i\in\bP(V) \mathrel{:} \prod p_i\neq 0\bigr\}$.
The dense orbit
$\mathcal O\bigl(\bigl[{\textstyle\sum v_i}\bigr]\bigr)$ in $X_{V,\cB}$
coincides with the intersection $X_{V, \cB} \cap \bT^{n-1}$. 
Observe that since $\dim X_{V,\cB}$ is equal to $\dim \mathcal
O\bigl(\bigl[\sum v_i\bigr]\bigr)$, it follows that  $\dim X_{V,\cB}$
is maximal among the dimensions
of the toric subvarieties of $\bP(V)$ --- i.e.\ those of the form $\overline{\mathcal
O([v])}$ for some $[v]\in \bP(V)$.

Based on the decomposition (\ref{eq:decomp}),  in
\cite[Proposition II.5.1.5]{kn:Gel}  it is proved
that any projective toric variety in an equivariant embedding is of type
$X_{V,\cB}$ for some $T$-module $V$ and a basis of eigenvectors
$\cB=\{v_1,\dots, v_n\}$ of $V$, in the
following sense. Let $U$ be a $T$-module and $Y\subset \bP(U)$ a toric
subvariety; then 
there exists $A=\{\lambda_1,\dots,\lambda_n\}\subset \cX(T)$ (with
possible repetitions) and a
$T$-equivariant linear injection $f: W:=\bigoplus_{i=1}^n
\Bbbk w _i \hookrightarrow U$, $t\cdot w_i=\lambda_i(t) w_i$, such that the induced
equivariant morphism $\widehat{f}:\bP(W)\hookrightarrow \bP(U)$ gives
an isomorphism $X_{W,\cB}\cong Y$.
Moreover, let $W'$ be another $T$-module, $\cB'=\{w'_1,\dots,w'_n\}\subset W'$  a basis of
eigenvectors of $W'$ such that  $t\cdot w'_i=\lambda_i(t)w'_i$, and consider $f\in
\Hom_T(W,W')$, given by
$f(w_i)=w'_i$. Clearly, $f$ is an isomorphism of $T$-modules, and its induced
morphism $\widehat {f}:\bP(W)\to \bP(W')$ is an isomorphism such that
$\widehat {f}(X_{W,\cB})=X_{W',\cB'}$.

In view of the preceding remark, the following notation makes sense.

\begin{definition} \label{def:vartor2}
The  projective toric
variety $X_A$ associated to the configuration of weights $A$ is
defined as
\[
X_A = X_{V,\mathcal B}= \overline{\cO\bigl(\bigl[\sum
v_i\bigr]\bigr)}\subset \bP(V), 
\]
where $V$ is a $T$-module with $A$ as associated configuration of weights.
\end{definition}

We can make a series of reductions on $A$ and $T$, as in \cite{kn:DFS05}.
First, the following easy lemma allows to reduce our problem to the
case of a faithful representation. 

\begin{lemma} \label{lem:faith}
Given a $T$-module $V$ of finite dimension and $A$ the associated
configuration of weights, consider the torus $T'= \Hom_\bZ (\langle A
\rangle_\bZ,\bk^*)$, where $\langle A\rangle_\bZ\subset \cX(T)$
denotes the $\bZ$-submodule generated by $A$.
The representation of $T$ in $\GL(V)$ induces a faithful representation $T' \to \GL(V)$ 
which has the same set theoretical orbits in $V$.
\qed
\end{lemma}

We can then replace $T$ by the torus $T'$. It is easy to show that
this is equivalent to the fact that $\langle A \rangle_\bZ=\cX(T)$,
which we will assume from now on without loss of generality.

\smallskip

Next, we  enlarge the torus without affecting the action on $\bP(V)$;
this will allow us to easily translate affine relations
to linear relations on the configuration of weights. 
If we let the algebraic torus $\Bbbk^*\times T$
act on $V$ by
$(t_0,t)\cdot v=t_0(t\cdot v)$, then the actions $T\times \bP(V)\to \bP(V)$ and
$(\Bbbk^*\times T)\times \bP(V)\to \bP(V)$ have  the same set theoretical
orbits. More in general, let $\lambda\in \cX(T)$ and
$A'=\{\lambda+ \lambda_1,\dots, \lambda+\lambda_n\}$. Consider the
$T$-action on $V$ given by $t\cdot_\lambda v_i=(\lambda+\lambda_i)(t)v_i$.
The actions $\cdot$ and $\cdot_\lambda$ coincide on
$\bP(V)$, and the corresponding variety $X_{A'}$ coincides with $X_A$.
Hence, we can assume that there is a splitting of
$T=\Bbbk^*\times S$ in such a way that $(t_0,s)\cdot v=t_0(s\cdot v)$ for all $v\in
V$,  $t_0\in \Bbbk^*$ and $s\in S$.

In fact, the previous  reductions  are comprised in the following more general
setting:

\begin{lemma}[{\cite[Proposition II.5.1.2]{kn:Gel}}]
\label{step:affine}
Consider $T,T'$ two tori and two
finite configurations of $n$ weights $A = \{ \lambda_1, \dots, \lambda_n\} \subset \cX(T)$, $A' =
\{ \lambda'_1, \dots, \lambda'_n\}\subset \cX(T')$.  
Assume that there exists a $\bQ$-affine transformation
$\psi: \cX(T)\otimes \bQ \to \cX(T') \otimes \bQ$ such that $\psi(\lambda_i) = \lambda'_i$
for all $i=1,\dots, n$. Then $X_A = X_{A'}$.
\qed
\end{lemma}

\begin{remark}
\label{rem:variado}
$(1)$ The dimension of the projective toric variety $X_A$ equals the
dimension of the affine span of $A$. 

\noindent $(2)$ Note that if $A=\{\lambda_1,\dots, \lambda_n\}$ is
contained in a hyperplane off the origin, then 
$X_A = \bP(V)$ precisely when $\dim T = n$ and the elements in $A$ are
a basis of $\cX(T)$.

\noindent (3) If we denote by $d$ the dimension of the affine 
span of 
$A$, then $X_A$ is a hypersurface if  and only if $n=d+2$. In this
situation, either $A$ coincides with the set of vertices of its convex
hull $\operatorname{Conv}(A)\subset \cX(T)\otimes \mathbb R$, or
$\operatorname{Conv}(A)$ contains 
only one element $\lambda \in A$ in its relative interior, and
$A\setminus \{\lambda\}$  is the set of vertices.
\end{remark}

We end this paragraph by recalling some basic facts about the
geometric structure of a toric variety $X_A$. 

\begin{lemma}[{\cite{kn:cls}}]
\label{lem:affineopen}
Let $A=\{\lambda_1,\dots, \lambda_n\}\subset \cX(T)$ be a
configuration, where $\{\lambda_1,\dots, \lambda_s\}$ is the set of
vertices of $\operatorname{Conv}(A)$. Set
$X_i=\operatorname{Spec}\bigl( \Bbbk\bigl[\mathbb Z^+ \langle(\lambda_j-\lambda_i)
\mathrel{:} \lambda_j\in A\rangle\bigr]\bigr)$,
$i=1,\dots, s$. Then $X_i$ is 
an affine toric $T$-variety, and there
exist $T$-equivariant  open immersions $\varphi_i:X_i\hookrightarrow
X_A$, in such 
a way that
\[
X_A=\cup_{i=1}^s \varphi_i(X_i)=\cup_{i=1}^s\operatorname{Spec}\bigl(
\Bbbk\bigl[\mathbb Z^+ \langle 
\lambda_j-\lambda_i\mathrel{:} \lambda_j\in A\rangle\bigr]\bigr).
\]

In particular, $X_A$ is a normal variety if and only if $\mathbb Z^+ \langle
\lambda_j-\lambda_i\mathrel{:} \lambda_j\in A\rangle=\bigl(\mathbb R^+ \langle
\lambda_j-\lambda_i\mathrel{:}\lambda_j\in A\rangle\bigr)\cap \cX(T)$ for
all $i=1,\dots ,s $.

Moreover, $X_A$ is a smooth variety if for all $i=1,\dots ,s$, 
there are exactly $\dim X_A$ edges of $\operatorname{Conv}(A)$ from $\lambda_i$,
and  the subset  $\{\lambda_{j_h} - \lambda_i\, : \,h =1 , \dots, \dim
X_A\}$ is  a
basis of $\cX(T)$, where $\lambda_{j_h}$ is the ``first'' point on an edge
from $\lambda_i$.
\qed 
\end{lemma}
\begin{proof}
See for example \cite[Appendix to Chapter 3]{kn:cls}.
\end{proof}

\subsection{Configurations in lattices, pyramids and projective joins} 

Let $M'$ be a lattice of
rank $d-1$. We let $M =\bZ \times M'$ and consider the $\bk$-vector space
$M_{\bk} = M\otimes_{\bZ}\Bbbk$.
Recall that given a basis $\{\mu_1,\dots ,\mu_d\}$ of $M$, we can identify $M$ with $\bZ^d$ and
$M_\bk$ with $\bk^d$.

\begin{definition} \label{def:regconf}
A {\em lattice configuration}\/ $A =\{\lambda_1, \dots,
\lambda_n\} \subset M$ is a finite sequence of lattice points.
We say that a configuration $A$ is  {\em regular}\/ if it is contained
in a hyperplane off the origin. 
\end{definition}

\begin{remark}
\label{exa:regularconf}
Let $T$ be an algebraic torus, and let
$A=\{\lambda_1,\dots, \lambda_n\} 
\subset \cX(T)$ be a  configuration of weights. Then the following are equivalent:

\noindent $(1)$  the configuration $A$ is regular;

\noindent $(2)$ up to affine isomorphism, $A$ has the form 
$\lambda_i = (1, \lambda'_i)$ for all $i = 1, \dots, n$;

\noindent (3) there exists a splitting
$T=\Bbbk^*\times S$, such that under the identification
$\cX(T)=\bZ\times \cX(S)$, the weights of $A$ are of the form
$\lambda_i=(1,\lambda' _i)$, $i=1,\dots,n$. See also the reductions
made before Lemma \ref{step:affine}.
\end{remark}

\begin{definition}
We denote by ${\mathcal R}_A\subset \bZ^n$ the lattice of affine relations
among the elements of $A$, that is  $(a_1,\dots,a_n)$
belongs to ${\mathcal R}_A$ if and only if $\sum_i a_i\lambda_i=0$ and
$\sum_i a_i=0$.  
\end{definition}

If  $A$ is regular, then  ${\mathcal R}_A$ coincides with the lattice
of linear relations  among the elements of $A$. Note that these
(affine or linear) relations 
among the elements of $A$ can be
identified with the  affine relations  
among the elements of the configuration $\{\lambda'_1, \dots,
\lambda'_n\} \subset M'$. 
Thus, given any configuration $\{\lambda'_1, \dots,
\lambda'_n\}\subset M'$, we can embed it 
in $M = \bZ \times M'$ via $\lambda' \mapsto (1, \lambda')$ so that
affine dependencies are 
translated to linear dependencies. In fact, the map $\lambda'\mapsto (1,
\lambda')$ is an injective affine 
linear map. More in 
general, we have the following definition.

\begin{definition}
We say that two configurations $A_i \subset \mathcal X(T_i), i =
1,2$, are \emph{affinely 
equivalent} if there exists an affine linear map $\varphi: \mathcal X(T_1)
\otimes \bR \to 
\cX(T_2) \otimes \bR$ (defined over $\bQ$) such that $\varphi$ sends
$A_1$ bijectively to $A_2$ (in particular, 
$\varphi$ defines an injective map from the affine span of $A_1$ to
the affine span of $A_2$). 
\end{definition}

So, if 
$A_1$ and $A_2$ are affinely equivalent, they have the same cardinal
and moreover,  
${\mathcal R}_{A_1} = {\mathcal R}_{A_2}$.  Any property of a
configuration $A$ shared by all its 
affinely equivalent configurations is called an {\em affine invariant}
of $A$.  In this terminology, Lemma \ref{step:affine} asserts that the
projective toric variety $X_A$ is an affine invariant of the
configuration $A$.

\begin{definition} \label{def:pyramid}
We say that  $A=\{\lambda_1,\dots ,\lambda_n\} \subset M$
is a {\em pyramid} \/ (or a {\em pyramidal configuration}\/) if there
exists an affine hyperplane $H$ 
such that $\# \{ i \, / \, \lambda_i \notin H\} = 1$, i.e. if all points in $A$ but one
lie in $H$, or equivalently, if there exist an index $i_0 \in \{1, \dots, n\}$ and
an affine linear function $\ell: M_\bk \to \bk$ such that $\ell(\lambda_i) =0$ for all $i \not= i_0$ and
$\ell(\lambda_{i_0}) =1$.
\par
More precisely, we say that $A$ is a {\em $k$-pyramidal configuration}\/ if, after
reordering, there exists a splitting of the lattice as a direct sum of lattices
 $M=M_1\oplus  M_2$, with
$A_1=\{\lambda_1,\dots ,\lambda_r\}$ a basis of $M_1$ and
$A_2=\{\lambda_{r+1},\dots ,\lambda_n\}\subset M_2$, with $A_2$ not
a pyramidal configuration of $M_2\otimes_{\bZ} \Bbbk$. In particular,
the $0$-pyramidal configurations are the \emph{non pyramidal
configurations}.
\end{definition}

\begin{remark} \label{obsnucleo}
Being a pyramid is clearly an affine invariant of a configuration.
 It is straightforward to check that $A$ is a non pyramidal
configuration if and only if 
there exists a relation $(p_1, \dots, p_n) \in \mathcal R_A$ with
$\prod_i p_i \neq 0$, i.e. 
if $\mathcal R_A$ is not contained in a coordinate hyperplane. 
\end{remark}

\begin{definition} \label{def:join}
Let $V_1, V_2$ two $\bk$-vector spaces of respective dimensions $h_1+1, h_2+1$ and 
$X \subset \bP(V_1), Y \subset \bP(V_2)$
two projective varieties. Recall that the  {\em join}\/ of  $X$
and $Y$ is the projective variety
 \[
\opJ_{h_1,h_2}(X,Y)= \overline{\bigl\{[x:y]:\,[x]
\in X, [y] \in Y \bigr\}} \subset \bP(V_1 \times V_2),
\]
that is, the cone over the join $\opJ_{h_1,h_2}(X,Y)$ is the
product of the cones over $X$ and $Y$.
We set
\[
\opJ_{h_1,h_2}(\emptyset,Y)=
\bigl\{[\underbrace{0:\dots:0}_{h_1+1}:y] \in
\bP(V_1\times V_2),\;[y] \in Y\bigr\} \subset
\bP(V_1 \times V_2).
\]
We define analogously $\opJ_{h_1,h_2}(X,\emptyset)$.
\par
We will denote $\bP^h = \bP(\bk^{h+1})$.
Observe that for any $Y \subset \bP^{h_2}$,  $Y\cong \opJ_{h_1,h_2}(\emptyset,Y)\subset
\bP^{h_1+h_2+1}$ for any $h_1 \in \bN$. 
If $X$ and $Y$ are non empty, then  $\dim
\opJ_{h_1,h_2}(X,Y)=\dim X+\dim Y+1$.

Remark that given $X_i\subset \bP(V_i)$, $\dim V_i=h_i+1$, $i=1,2,3$,
then 
\[
\opJ_{h_1+h_2+1, h_3}\bigl(\opJ_{h_1,h_2}(X_1,X_2),X_3\bigr)=
\opJ_{h_1,h_2+h_3+1}\bigl(X_1,\opJ_{h_2,h_3}(X_2,X_3)\bigr)\subset \bP(V_1\times
V_2\times V_3).
\]
We will denote this variety by $\opJ_{h_1,h_2,
  h_3}(X_1,X_2,X_3)$.
\end{definition}

Given two projective toric varieties $X_{A_1}$ and $X_{A_2}$, then
their join  is also a toric variety:

\begin{remark}
\label{obse:parajoin}
$(1)$ Let $T=S_1\times S_2$ be a splitting of $T$ as a product of tori, and
$A_1=\{\lambda_1,\dots,\lambda_k\}\subset \cX(S_1)$,
$A_2=\{\lambda_{k+1},\dots,\lambda_{n}\} \subset  \cX(S_2)$ two regular configurations.

Let $V_1=\bigoplus_{i=1}^{k} \Bbbk v_i$, $s_1\cdot
v_i=\lambda_i(s_1)v_i$ for all $s_1\in S_1$, and
$V_2=\bigoplus_{i=k+1}^{n} \Bbbk v_i$, $s_2\cdot
v_i=\lambda_i(s_2)v_i$ for all $s_2\in S_2$.
Then $V=V_1\times V_2$ is a $T$-module for the product action $(s_1,s_2)\cdot
(w_1,w_2)=(s_1\cdot w_1, s_2\cdot w_2)$. Moreover, $V$ decomposes in
simple submodules as $V=\bigoplus_{i=1}^k \Bbbk (v_i,0)\oplus
\bigoplus_{i=k+1}^n \Bbbk(0,v_i)$.

Consider the
$S_i$-toric varieties $X_{A_i}\subset \bP(V_i)$ ($i=1,2$) and let $A=A_1\times\{0\}\cup
\{0\}\times A_2\subset \cX(S_1)\times \cX(S_2)=\cX(T)$. 
The projective toric variety associated to $A$ is then the join
$X_A=\opJ_{k-1,n-k-1}(X_{A_1}, X_{A_2})$.

\noindent $(2)$ In the particular case when $A\subset M=M_1\oplus M_2$
is a $k$-pyramidal configuration with $A_1\subset  
M_1, A_2\subset M_2$ as in Definition~\ref{def:pyramid},
let $ S_1 = \Hom_\bZ (M_1, \bk^*)$, $S_2 = \Hom_\bZ (M_2, \bk^*)$, $T
= S_1 \times S_2$ and $V$ 
as above. We then have that $X_A = \opJ_{k-1, n-k-1}(\bP(V_1),
X_{A_2})$. That is, $X_A$ is the \emph{cone over $X_{A_2}$ with vertex
  $\bP(V_1)$}. 
\end{remark}

Next, we describe the toric varieties associated  to configurations 
with repeated weights. Recall that a projective variety is called {\em non
degenerate} if it is not contained in a proper linear subspace.

\begin{lemma} 
\label{lema:splittingjoin} 
Let $A=\{\lambda_1,\dots,\lambda_1,\dots,\lambda_h,\dots,
\lambda_h\}\subset \cX(T)$ be a configuration of $n$ weights, with
$\lambda_i$ appearing $k_i+1$ times  and $\lambda_i\neq\lambda_j$ if
$i\neq j$. If we set  $k= \sum_i k_i=n-h$, then the smallest linear
subspace that contains $X_A$ has codimension  $k$.

In particular, $X_A$ is a non degenerate variety if and only if the
configuration $A$ has no repeated elements.
\end{lemma}

\begin{proof}
Let $\cB=\{v_{1,1},\dots, v_{1,k_1+1},\dots, v_{h,1},\dots, 
v_{h,k_h+1}\}$ be a basis of associated eigenvectors of $V$, with $t\cdot
v_{i,j_i}=\lambda_i(t)v_{i,j_i}$ for all $i=1,\dots ,h$, $j_i=1,\dots, k_i+1$.
 Consider a  hyperplane $\Pi\subset \bP\Bigl(
\bigoplus_{i=1}^h\bigl(\bigoplus_{j_i=1}^{k_i+1}
\Bbbk v_{i,j_i}\bigr)\Bigr)$ of equation
\[
\sum_{i,j_1,\dots, j_h}
c_{i,j_i}x_{i,j_i}=0,
\]
where $x_{i,j_i}$ are the coordinates in the basis $\cB$. 
Then
$X_A\subset \Pi$ if and only if $\bigl[t\cdot \sum
v_{i,j_i}\bigr]\subset \Pi$ for all $t \in T$. As 
$\bigl[t\cdot \sum v_{i,j_i}\bigr]=\bigl[\sum \lambda_i(t) v_{i,j_i}\bigr] \in \Pi$, this
 is equivalent to the equalities
\[
\sum_{i=1}^h \sum_{j_i=1}^{k_i+1} c_{i,j_i}\lambda_i(t)= 0, \, \, t \in T.
\]

Since $\{\lambda_1,\dots,\lambda_h\}$ are different weights, we
deduce that $\sum_{j_i=1}^{k_i+1}c_{i,j_i}=0$ for all $i=1,\dots,
h$.  It follows that the maximum codimension of a subspace that
contains $X_A$ is $\sum_{i=1}^h k_i=k$.

On the other hand, clearly
\[
X_A\subset H \, = \,
\bigl\{{\textstyle \sum_{i=1}^h
x_i \sum_{j_i=1}^{k_i+1} v_{i,j_i}}
\mathrel{:} x_i\in \Bbbk
\bigr\},
\]
where the subspace $H\subset \bP(V)$ has codimension $k$.
\end{proof}

\begin{lemma}
\label{lem:otromas}
Let $A=\{\lambda_1,\dots,\lambda_1,\dots,\lambda_h,\dots, 
\lambda_h\}\subset \cX(T)$ be a configuration of $n$ weights, with
$\lambda_i$ appearing $k_i+1$ times   and
$\lambda_i\neq\lambda_j$ if $i\neq j$. Set  $k= \sum_i k_i=n-h$ and let
\[
V=\bigoplus_{i=1}^h\bigr({\textstyle \bigoplus_{j_i=1}^{k_i+1} \Bbbk
v_{i,j_i}}\bigr)= \Bigr(\bigoplus_{i=1}^h\bigr({\textstyle
\bigoplus_{j_i=1}^{k_i} \Bbbk 
v_{i,j_i}}\bigr)\Bigr) \oplus \Bigr({\textstyle \bigoplus_{i=1}^{h} \Bbbk
v_{i,j_{k_i +1}} }\Bigr), 
\]
with $t\cdot
v_{i,j_i}=\lambda_i(t)v_{i,j_i}$ for all $t\in T$, $i=1,\dots h$,
$j_i=1,\dots, k_i+1$.

Let $C=\{\lambda_1,\dots,\lambda_h\}$  and  consider  $X_C\subset \bP
\bigl(\bigoplus_{i=1}^{h} \Bbbk
v_{i,j_{k_i +1}} \bigr)$. Then
$X_A$ is isomorphic to the cone $\opJ_{k-1,h-1}  (\emptyset, X_C)$ over the non
degenerate projective toric variety $X_C$. 
\end{lemma}

\begin{proof}
Let $f:V \to V$ the linear isomorphism defined by
\[ 
\begin{split}
f\bigl( (x_{i,j_i})_{i=1,\dots,h, j_i=1, \dots, k_i}, (x_{i,
j_{k_i+1}})_{i=1,\dots,h}\bigr) \, = \\ 
 \bigl( (x_{i,j_i}- x_{i, j_{k_i+1}})_{i=1,\dots,h, j_i=1, \dots,
k_i}, (x_{i,j_{ k_i+1}})_{i=1,\dots,h}\bigr). 
\end{split}
\]
The associated projective map clearly sends
$X_A$ to the join $\opJ_{k-1,h-1}  (\emptyset, X_C)$.
\end{proof}

In Proposition \ref{prop:jointoric} below we combine Remark \ref{obse:parajoin}
and Lemmas \ref{lema:splittingjoin}  and 
\ref{lem:otromas}, in order to describe a projective toric variety as a cone
over a non degenerate projective \emph{toric} variety that is not a
cone (that is, the associated configuration is non pyramidal).

\begin{remark}
\label{rem:Xnd}
Let $X\subset \bP^{n-1}$ be a \emph{non linear} irreducible projective
variety.
Let $H\subset \bP^{n-1}$ be the minimal linear subspace containing
$X$, and let $k$ be the codimension of $H$. Then $H\cong
\bP^{n-k-1}$ and if $X'$ denotes the variety 
$X$ as a subvariety of $H$, then $X=\opJ_{k-1,n-k-1}(\emptyset,
X')$. Since $X'$ is non degenerate, it follows that there exists
$Y\subset \bP^{m-1}$  such that
$X'=\opJ_{h-1,m-1}\bigl(\bP^{h-1},Y\bigr)$, where $n-k-1= h +m
-1$. Hence, we have an identification
\[
X=\opJ_{k-1,h-1,m-1}\bigr(\emptyset, \bP^{h-1},
Y\bigr).
\]

In particular, $\dim X=h+ \dim Y$.

Observe  that $Y\subset \bP^{m-1}$ is a non degenerate
subvariety. Moreover, we can assume that $Y$ is not a cone. 
In this
case, we will denote $X_{\operatorname{nd}}=Y$. If moreover $X$ is an
equivariantly embedded toric variety, then we can choose $X_{\nd}$ as
$X_{C_2}$ in the following proposition.

When $X$ is linear, $X=H$, $m=1$  and $Y$ is empty.
\end{remark}

\begin{proposition}
\label{prop:jointoric}
Let $A=\{\lambda_1,\dots,\lambda_1,\dots,\lambda_h,\dots, 
\lambda_h\}\subset \cX(T)$ be a configuration of $n$ weights, with
$\lambda_i$ appearing $k_i+1$ times  and
$\lambda_i\neq\lambda_j$ if $i\neq j$. Set  $k= \sum_i k_i=n-h$ and assume that
$C=\{\lambda_1,\dots,\lambda_h\}$   is $r$-pyramidal. Then there
exists a splitting  $T=S_1\times S_2$  such that, after 
reordering of the elements in $C$, it holds that $C=C_1\cup C_2$, where
$C_1=\{\lambda_1,\dots,\lambda_r\}$  is a  basis of $\cX(S_1)$ and
$C_2=\{\lambda_{r+1},\dots,
\lambda_{h}\}\subset  \cX(S_2)$ is  a non pyramidal configuration, as in 
Definition~\ref{def:pyramid}. Moreover, we have that 
\[
X_A=\opJ_{k-1,r-1,h-r-1}\bigr(\emptyset, \bP^{r-1},
X_{C_2}\bigr).
\]
In the special case when $X_A$ is linear, $C_2$ is empty.
\end{proposition}
\begin{proof}
 We set $V=\bigoplus_{i=1}^h\bigr(\bigoplus_{j_i=1}^{k_i+1} \Bbbk
v_{i,j_i}\bigr)$, with $t\cdot
v_{i,j_i}=\lambda_i(t)v_{i,j_i}$ for all $t\in T$, $i=1,\dots h$,
$j_i=1,\dots, k_i+1$, and $w_i = v_{i,k_i+1}$.

Assume that  $C$ is a $r$-pyramidal
configuration, and let  $X_C\subset
\bP\bigr(\bigoplus_{i=1}^h \Bbbk w_i\bigl)$. Then,  there exists a splitting
$T=S_1\times S_2$ such that, after reordering of $C$,
$C_1=\{\lambda_1,\dots,\lambda_r\}$  is a  basis of $\cX(S_1)$ and
$C_2=\{\lambda_{r+1},\dots,
\lambda_{h}\}\subset  \cX(S_2)$ is  a non pyramidal configuration. Hence,
\[
X_C=\opJ_{r-1,h-r-1}
(X_{\Id_r},X_{C_2})= \opJ_{r-1,h-r-1}
\bigl(\bP\bigl(\oplus_{i=1}^r
\Bbbk v_{i,k_i+1}\bigr),X_{C_2}\bigr).
\]

By Lemma \ref{lem:otromas}, we  can assume that
$X_A=\opJ_{k-1,h-1}
(\emptyset ,X_{C})$, and so  
\[
\begin{split}
X_A= &\ \opJ_{k-1,h-1}\bigl(\emptyset,\opJ_{r-1,h-r-1}
(X_{\Id_r},X_{C_2})\bigr)= \\
&\ 
 \opJ_{k-1,h-1}\bigl(\emptyset, \opJ_{r-1,h-r-1}
\bigl(\bP\bigl(\oplus_{i=1}^r
\Bbbk v_{i,k_i+1}\bigr),X_{C_2}\bigr)\bigr)= \\
&\ \opJ_{k-1,r-1,h-r-1}\bigl(\emptyset, 
\bP\bigl(\oplus_{i=1}^r
\Bbbk v_{i,k_i+1}\bigr),X_{C_2}\bigr),
\end{split}
\]
as claimed.
\end{proof}

\subsection{Dual of a projective toric variety}

We recall the classical notion of the dual variety of a projective variety.

\begin{definition} \label{def:dual}
Let $V$ be a $\Bbbk$-vector space of finite dimension and denote by
$V^\vee$ its dual $\bk$-vector space. Let $X\subset
\bP(V)$ be an irreducible  projective variety. The
{\em dual variety}\/ of $X$ is defined as the closure of
 the hyperplanes intersecting the regular part $X_{reg}$ of
 $X$ non transversally:
\[
X^*=\overline{\bigl\{ [f]\in \bP(V^\vee) \mathrel{:} \exists x\in
X_{reg}\,,\ f|_{T_xX}\equiv 0\bigr\}}\subset \bP(V^\vee).
\]
As usual, $T_xX$ denotes the embedded tangent space of $X$ at $x \in X_{reg}$.

Note that $\bP(V)^* = \emptyset$. We set by convention,
$\emptyset^*=\bP(V^\vee)$. 
\end{definition}

Self-duality is not an intrinsic property, it
depends on the projective embedding. It can be proved that $X^*$ is an
irreducible projective variety and that $(X^*)^*=X$ 
(see for example \cite{kn:Gel}). 

For a generic variety $X\subset
\bP(V)$,  $\codim X^*=1$. 
If $\codim X^*\neq 1$, it is said that $X$ has {\em defect}\/   $\codim
X^*-1$.

\begin{definition} \label{def:selfdual}
An irreducible projective variety $X \subset \bP(V)$ is called 
{\em self-dual}\/ if $X$ is isomorphic to $X^*$ as embedded projective
varieties, that is if there exists a (necessarily linear)
isomorphism 
$\varphi:\bP(V)\to \bP(V^\vee)$ such that $\varphi(X)=X^*$.
\end{definition}

A self-dual projective variety $X \subset \bP^{n-1}$ of dimension $d -1< n-1$
(i.e., which is not a hypersurface) has  positive defect
$n-d-1$. The defect of the whole projective space
$\bP^{n-1}$ is $n-1$.

\begin{remark}\label{rk:identesp-espdual}
\noindent Recall that given a basis $\cB=\{v_1,\dots,v_n\}$ of
$V$,  we can identify
$\bP(V)$ with $\bP(V^\vee)$ by means of $v_i \leftrightsquigarrow v_i^\vee$,
 where $\{v_1^\vee,\dots,v_n^\vee\}$ is the dual basis of $\cB$.
Then, via the choice of a basis of $V$, we can look at the dual variety
inside the same projective space. Self-duality
can be reformulated as follows: $X\subset \bP(V)$ is self-dual if there
exists  $\varphi\in \Aut\bigl(\bP(V)\bigr)$ such  that
$\varphi(X)=X^*$.
\end{remark}

Let $V$ be a $T$-module of finite dimension $n$ over a $d$-dimensional torus $T$ 
and let $A$ be the associated configuration of weights.
In view of the considerations of the preceding subsections, we assume from now
on and without loss
of generality, that $A=\{\lambda_1,\dots,\lambda_n\}\subset \cX(T)$ is a regular 
configuration, possibly with repeated elements, such
that $\langle A \rangle_\bZ=\cX(T)$. 

The regularity of $A$ implies in particular the existence 
of a splitting $T = \bk^* \times S$ as in 
Remark  \ref{exa:regularconf}. Then, $X_A$ is a  $(d-1)$-dimensional subvariety
of the $(n-1)$-dimensional projective space $\bP(V)$ and the lattice
${\mathcal R}_A$ has rank $n-d$.

The dual variety $X_A^*$ has the following interpretation.
For $[\xi]\in\bP(V^\vee)$, let  $f_{\xi}\in \Bbbk[T]$, $f_\xi(t)=\xi(t\cdot\sum
v_i)\in\Bbbk[T]$.
Then $X_A^*$ is obtained as the closure of the set of those
$[\xi]\in\bP(V^\vee)$ such that
there exists $t\in T$ with $f_\xi(t)=\frac{\partial f_\xi}{\partial
t_i}(t)=0$ for all $i=1,\dots ,n$.

\[
X_{A}^{*}=\overline{\left\{\xi \in
\bP(V^\vee) \mathrel{:} \exists t\in T\,,
f_{A}(t)=\frac{\partial f_\xi}{\partial
t_1}(t)=\frac{\partial f_\xi}{\partial t_2}(t)=\dots=\frac{\partial
f_\xi}{\partial t_d}(t)=0\right\}}.
\]

In \cite{kn:DFS05} a rational parameterization of the dual variety $X_A^*$ was
obtained. We adapt this result to our notations.
As before, $\cB=\{v_1,\dots, v_n\}$ is a basis of
eigenvectors, $t\cdot v_i=\lambda_i(t)v_i$, and $\cB_A=\{u_1,\dots,u_{n-d}\}$ is a basis of 
$\mathcal R_A$.  We denote by ${\mathcal R}_{A, \bk}$ the $(n-d)$-dimensional $\bk$-vector space 
${\mathcal R}_A \otimes_\bZ \bk$ and we identify
 $\bP(V)$ with $\bP(V^\vee)$  by means of the chosen basis
$\mathcal B$ of eigenvectors (and its dual basis) as in Remark~\ref{rk:identesp-espdual}.

\begin{proposition}[{\cite[Proposition~4.1]{kn:DFS05}}]
\label{thm:dualvarparam} 
Let $T=\Bbbk^*\times S,V,A,\cB, \cB_A$ as before.
Then the mapping $\bP({\mathcal R}_{A, \bk})\times
S\to \bP(V)$ defined by 
\[\bigl([a_1:\dots : a_{n}],s) \mapsto  s\cdot\bigl[\sum
a_iv_i\bigr] 
\]
has image dense in $X_A^*$. That is, the morphism
\[
(\Bbbk^*)^{n-d}\times
T\to \bP(V)\,,  \quad (c,t) \mapsto  t\cdot\bigl[\sum c_iu_i\bigr]
\]
is a rational parameterization of $X_A^*$, and 

\[
X_{A}^{*}=\overline{\bigcup \limits_{p \in \bP
(\mathcal R_{A,\bk})} \mathcal{O}\bigl(p \bigr)}=\overline{T\cdot
\bP(\mathcal R_{A,\bk})}.\vspace*{-0.5cm}
\]
\qed
\end{proposition}

This last equality, which expresses the dual variety as the closure of the union of the torus
orbits of all the classes in the vector space of relations of the configuration $A$,
 is the starting point of our classification of self-dual projective
toric varieties, which we describe in the sequel.

\section{Characterization of self-duality in terms of orbits} \label{sec:orbits}

Let $T$ be a torus of dimension $d$ and  $V$ a rational $T$-module of
dimension $n$ with associated configuration of weights  $A =
\{\lambda_1,\dots, \lambda_n\}$. We assume that  
 $\langle A \rangle_\bZ=\cX(T)$ and  keep the notations of the
preceding section. Given $p=\bigl[\sum p_iv_i\bigr]\in \bT^{n-1}$,  
we denote by $m_p\bigl([\sum
x_iv_i]\bigr)=\bigl[\sum p_ix_iv_i\bigr]$ the diagonal linear
isomorphism defined by $p$.

\subsection{Non pyramidal configurations} 

In this subsection we characterize self-dual projective toric varieties
associated to a configuration of weights $A$ which define a {\em non pyramidal} \/
configuration, in terms of the orbits of the torus action.

Note that the whole projective space $\bP(V)$ can be seen as a toric 
projective variety associated to a $\dim V$-pyramidal configuration 
and its dual variety is empty. But we now show
that for non pyramidal configurations the dimension of the dual variety
$X_A^*$ cannot be smaller than the dimension of the toric variety $X_A$.
This result has been proved by Zak \cite{kn:zak} for any non degenerate smooth
projective variety.

\begin{lemma}
\label{lema:nonpiramidaldim}
If $A$ is a non pyramidal configuration, then $\dim X_A^*\geq \dim
X_A$.
\end{lemma}

\begin{proof}
Indeed,  if $A$ is not a pyramidal
configuration, then by Remark \ref{obsnucleo} we know that there
exists $p=\sum p_iv_i\in \mathcal R_{A,\bk}$ such
that $p_i\neq 0$ for all $i=1,\dots, n$. Hence, if we identify
$\bP(V)$ with $\bP(V^\vee)$ by means of the dual basis, then
\[
X_A^*=\overline{T\cdot \bP(\mathcal R_{A,\bk})}\supset
\overline{\cO ([p])}= m_p\Bigl(\overline{\cO\bigl(
\bigl[{\textstyle \sum_iv_i}\bigr]\bigr)}\Bigr) = m_p(X_A).
\]
Since $p \in \bT^{n-1}$, we have that $\dim m_p (X_A) = \dim X_A$ and
the result follows.
\end{proof}

We identify $\bP(V)$ with $\bP(V^\vee)$  by means of the chosen basis
$\mathcal B$ of eigenvectors (and its dual basis) as in Remark~\ref{rk:identesp-espdual}.
The following is the main result of this subsection.

\begin{theorem}
\label{teoeqauto}
Let $A\subset \cX(T)$ be a non pyramidal configuration. 

The following assertions are equivalent.

\begin{enumerate}
\item[(1)] $X_A$ is a self-dual projective variety.

\item[(2)] There exists  $p_0 \in \bP
(\mathcal R_{A,\bk})\cap \bT^{n-1}$
such that  $\bP (\mathcal R_{A, \bk})\subset
\overline{\mathcal{O}(p_0)}$.

\item[(3)] There exists  $p_0 \in \bP
(\mathcal R_{A, \bk})\cap
\bT^{n-1}$ such that  $X_{A}^{*} = m_{p_0}( X_{A})$.

\item[(4)]  For all $q \in \bP (\mathcal R_{A, \bk})\cap \bT^{n-1}$, $\bP
(\mathcal R_{A, \bk})\subset \overline{\mathcal{O}(q)}$.

\item[(5)] For all $q \in \bP (\mathcal R_{A,k})\cap \bT^{n-1}$,
$X_{A}^{*} =m_q(X_{A})$.
\end{enumerate}
\end{theorem}

\begin{proof}
We prove $(1) \Rightarrow (5)$ and $(2) \Rightarrow (4)$, the rest
of the implications being trivial.

\noindent $(1) \Rightarrow (5)$: By Proposition \ref{thm:dualvarparam},
\[
\begin{split}
X_{A}^{*}= & \overline{\bigcup_{p \in \bP
(\mathcal R_{A, \bk})} \mathcal{O}(p)}\supset \overline{\bigcup_{p \in
\bP (\mathcal R_{A, \bk}) \cap \mathbb{T}^{n-1}}\mathcal{O}(p)} \supset
\overline{\mathcal{O}(q)}= m_{q}(X_A),
\end{split}
\]
for all  $q \in \bP (\mathcal R_{A,\bk}) \cap \bT^{n-1}$. Since $\dim X_A=\dim
X_A^*$, equality holds in the last equation.

$(2) \Rightarrow (4)$:  Let  $p_0 \in \bP (\mathcal R_{A,\bk})\cap \bT^{n-1}$
be such that  $\bP (\mathcal R_{A,\bk})\subset
\overline{\mathcal{O}(p_0)}$. If  $q\in  \bP (\mathcal R_{A,\bk})\cap \bT^{n-1}$, then
$q\in \overline{\mathcal{O}(p_0)}\cap \bT^{n-1}=\cO(p_0)$. Then,
${\mathcal O}(q) = {\mathcal O}(p_0)$ and the result follows.
\end{proof}

The equivalence between $(1)$ and $(5)$ in Theorem \ref{teoeqauto} implies that as 
soon as the dual of an equivariantly embedded projective toric variety of the form $X_A$
has the same dimension of the variety, there exists a linear isomorphism between them.

\begin{theorem}\label{thm:autdim}
Let $A\subset \cX(T)$ be a configuration of weights  which is non pyramidal.
 Then $X_A$ is self-dual if and only if $\dim X_A=\dim
X_A^*$.
\qed
\end{theorem}

This result is not true in general
for projective toric varieties not equivariantly embedded, even for rational planar
curves (for which the dual is again a curve, but not necessarily isomorphic).

\subsection{The general case} 

We now address the complete
characterization of  self-dual projective toric
varieties associated to an arbitrary configuration of weights $A\subset
\cX(T)$. We keep the notations of the preceding section. 

We  begin by recalling a well known result about duality of
projective varieties:

\begin{lemma}[{\cite[Theorem 1.23]{kn:tev05}}]
\label{lem:tevelev}
Let $X \subset \bP^n$ be  a non linear irreducible subvariety.

$(1)$ Assume that $X$ is contained in a hyperplane $H = \bP^{n-1}$. If
${X'}^*$ is the 
    dual variety of $X$, when we consider $X$ as a subvariety of
    $\bP^{n-1}$, then $X^*$ is the cone over ${X'}^*$ with vertex $p$
    corresponding to $H$. 

$(2)$ Conversely, if $X^*$ is a cone with vertex $p$, then $X$ is
contained in the  corresponding hyperplane $H$.
   
When $X$ is linear, $(X')^*$ is empty.    
\qed
\end{lemma}

As an immediate application of Lemma~\ref{lem:tevelev}, we have the following characterization of
self-dual equivariantly embedded projective toric hypersurfaces. Note that the only linear varieties which
are self-dual are the subspaces of dimension $k-1$ in $\bP^{2k-1}$. In particular,
the only hyperplanes which are self dual are points in $\bP^1$.

\begin{corollary}
\label{coro:noname}
Let $T$ be an algebraic torus and $A\subset \cX(T)$ a configuration such that
$X_A$ is a non linear hypersurface. Then
$X_A$ is self-dual if and only if $X_A$ is not a cone.
\end{corollary}

\begin{proof}
Assume that $X_A$ is a cone. Then by Lemma \ref{lem:tevelev}, it
follows that $X_A^*$ is contained in a hyperplane, hence $X_A$ is not
self-dual. 

If $X_A$ is not a cone, then $A$ is non pyramidal (see Remark \ref{obse:parajoin}), and it follows from
Lemma \ref{lema:nonpiramidaldim} that $\dim
X_A^* \geq \dim X_A$. If $\dim
X_A^* > \dim X_A$, then $X_A^*=\bP(V)$ and hence
$X_A=(X_A^*)^*=\emptyset$, which is a contradiction. It follows that $\dim
X_A^* = \dim X_A$ and hence  
Theorem \ref{thm:autdim} implies that $X_A$ is self-dual.  
\end{proof}

Applying Lemma \ref{lem:tevelev}, we can reduce the study of  duality
of  projective  varieties  to the study of non degenerate projective 
varieties that are not a  cone.

\begin{proposition}
\label{prop:arb}
Let $X\subset \bP^{n-1}$ be an  irreducible projective
variety. Let $k-1$ be the codimension of the minimal subspace of
$\bP^{n-1}$ containing $X$. Then, with the notations of Remark
\ref{rem:Xnd}, the following assertions hold:

\noindent $(1)$ If  $X=\opJ_{k-1,k-1,m-1}\bigr(\emptyset,
\bP^{k-1}, 
X_{\nd}\bigr)$, with $X_{\nd}\subset \bP^{m-1}$ self-dual, then $X$ is self-dual.

\noindent $(2)$  If $X$ is self-dual, then $\dim X_{\nd}=\dim
(X_{\nd})^*$, and $h=k$, that is 
\begin{equation}\label{eq:sd}
X=\opJ_{k-1,k-1,m-1}\bigr(\emptyset,
\bP^{k-1}, 
X_{\nd}\bigr),
\end{equation} 
\end{proposition}

\begin{proof}
Let $X=\opJ_{k-1,h-1,m-1}\bigr(\emptyset, \bP^{h-1},
X_{\nd}\bigr)$.  Applying recursively
Lemma \ref{lem:tevelev} (see Remark \ref{rem:Xnd}) we obtain that 
\[
\begin{split} 
X^*=&\ \opJ_{k-1,h-1,m-1}\bigr(\emptyset, \bP^{h-1},
X_{\nd}\bigr)^*=\opJ_{k-1,h+m-1}\bigr(\emptyset, \opJ_{h-1,m-1}\bigl(\bP^{h-1},
X_{\nd}\bigr)\bigr)^*= 
\\
& \ 
\opJ_{k-1,h+m-1}\bigr(\bP^{k-1}, \opJ_{h-1,m-1}\bigl(\bP^{h-1},
X_{\nd}\bigr)^*\bigr)= \\
& \  \opJ_{k-1,h+m-1}\bigr(\bP^{k-1}, \opJ_{h-1,m-1}\bigl(\emptyset,
X_{\nd}^*\bigr)\bigr)=\\
& \  \opJ_{k-1,h-1,m-1}\bigr(\bP^{k-1}, \emptyset, 
X_{\nd}^*\bigr)= \opJ_{h-1,k-1,m-1}\bigr( \emptyset, \bP^{k-1}, 
X_{\nd}^*\bigr)=\\
& \ 
\opJ_{h-1,k+m-1}\bigr(\emptyset, \opJ_{k-1,m-1}\bigl(\bP^{k-1},
X_{\nd}^*\bigr)\bigr).
\end{split}
\]

In particular, $\dim X^*= k+\dim X_{\nd}^*$, and the maximal subspace 
that contains $X^*$ has codimension $h$.

In order to prove $(1)$, assume  that $h=k$ and  $X_{\nd}$ is 
self-dual. Then $X^*= 
\opJ_{k-1,k-1,m-1}\bigr(\emptyset, 
\bP^{k-1}, 
X_{\nd}^*\bigr)$.  Since $X_{\nd}$ is self-dual, there exists an
isomorphism $\varphi: \bP^{m-1}\to  \bP^{m-1}$ such that
$\varphi(X_{\nd})=X_{\nd}^*$. It is clear that  $\varphi$ extends to an
isomorphism $\widetilde{\varphi}:\bP^{n-1}\to
\bP^{n-1}$ such that $\widetilde{\varphi}(X)=X^*$.

In order to prove $(2)$, 
assuming $X$ is self-dual and
writing $X$ as in Remark~\ref{rem:Xnd}, 
it follows that $h=k$, and hence $h+\dim X_{\nd}=\dim X= \dim X^*=k+\dim X_{\nd}$.\end{proof}

In our toric setting, Proposition \ref{prop:arb} can be improved, so
that we obtain a geometric
characterization of self-dual projective toric varieties.

\begin{theorem}\label{thm:autdimgral}
Let $A$ be an arbitrary lattice configuration. Then $X_A$ is self-dual if and
only if  $\dim X_A=\dim X_A^*$ and the smallest linear subspaces containing
$X_A$ and $X_A^*$ have the same (co)dimension.
\end{theorem}

\begin{proof} 
By Proposition \ref{prop:jointoric}, 
\[
X_A=\opJ_{k-1,h-1,m-1}\bigr(\emptyset, \bP^{h-1},
X_{C_2}\bigr),
\]
where $C_2\subset A$ is a non pyramidal configuration without
repeated weights. By Theorem \ref{thm:autdim}, $X_{C_2}\subset
\bP^{m-1}$ is self dual if and only if $\dim X_{C_2}=\dim
X_{C_2}^*$. The result follows now from Proposition \ref{prop:arb}.
\end{proof}

Combining Proposition \ref{prop:jointoric} and Theorem  
\ref{thm:autdimgral} we obtain the following explicit combinatorial
description 
of self-dual toric varieties.

\begin{theorem}
\label{teo:constaut+gde}
Let $A=\{\lambda_1,\dots ,\lambda_1,\lambda_2,\dots, \lambda_2,\dots,
\lambda_h,\dots,\lambda_h\}\subset \cX(T)$  be a configuration of $n$ weights with  each
$\lambda_i$ appearing $k_i+1$ times, $\lambda_i\neq \lambda_j$ if
$i\neq j$. Let
$C=\{\lambda_1,\dots ,\lambda_h\}$ be the associated configuration without repeated
weights. Then $X_A$ is self-dual if and only
if the following assertions hold.

$(1)$ $C$ is a $k$-pyramidal configuration, where $k=n-h=\sum k_i$.

$(2)$ There exists a splitting  $T=S_1\times S_2$  such that, after
reordering of the elements in $C$, it holds that $C=C_1\cup C_2$, where
$C_1=\{\lambda_1,\dots,\lambda_k\}$  is a  basis of $\cX(S_1)$ and
$C_2=\{\lambda_{k+1},\dots,
\lambda_{h}\}\subset  \cX(S_2)$ is  a non pyramidal configuration, as in 
Definition~\ref{def:pyramid}. Moreover, the $S_2$--toric projective variety
$X_{C_2}\subset \bP\bigl(\bigoplus_{i=k+1}^{h}  \Bbbk w_{i}\bigr)$,
$t\cdot w_i=\lambda_i(t)w_i$, is self-dual.
\qed
\end{theorem}

It follows from Theorem \ref{teo:constaut+gde} that if $X_A$ is a
self-dual toric variety with $A$ pyramidal, then there are
repeated weights in $A$. The converse of this statement does not hold.
In the next example we show a family of {\em non}-pyramidal
configurations $A$ with repetitions such that $X_A$ is self dual.

\begin{example}\label{ex:conversenotrue}
Let $C=\{c_1, \dots, c_s\}\subset \bZ^{n-1}$  be  any non pyramidal
configuration, such that $X_C$ is self-dual. Then, the configuration
$A = \{ e_1, e_1, (0, c_1), \dots, (0, c_s) \}\subset \bZ^n$ has
repeated weights, and $X_A$ is self-dual by Theorem
\ref{teo:constaut+gde}. It is straightfoward to check that $A$ is non-pyramidal.
Note that these configurations
become pyramidal when we avoid repetitions.
\end{example}

\section{Characterizations of self-duality in combinatorial terms}
\label{sec:comb}

In this section we will characterize self-duality of projective
toric varieties of type $X_A$ in combinatorial terms. 
We make explicit calculations for the algebraic torus $(\Bbbk^{*})^d$
acting on $\Bbbk^n$, in order to give an interpretation of
the conditions of Theorem \ref{teoeqauto} in terms of the configuration
$A$ and in terms if its {\it Gale
dual configuration\/}, whose definition we recall below. 

We refer the reader to 
\cite[Chapter 6]{kn:Zi} for an account of the basic combinatorial
notions we use in what follows.

\subsection{Explicit calculations for ${(\Bbbk^{*})}^d$ acting on
$\Bbbk^n$} \label{rk:ector}

Let $T=(\Bbbk^{*})^d$. We identify the lattice of characters 
$\mathcal X(T)$ with $\mathbb{Z}^d$. Thus, any character $\lambda\in
\mathcal{X}(T)$ is of the form $\lambda(t)=t^m$, 
where $m \in \bZ^d$ and $t^m=t_1^{m_1}\cdots t_d^{m_d}$.  We take the
canonical basis of $\Bbbk^n$ as the basis of eigenvectors of the action
of $T$. That is, if 
$A= \{\lambda_1,\dots,
\lambda_n\}\subset \bZ^d$, $T$ acts on $\Bbbk^n$
by $t\cdot (z_1,\dots,z_n)=(t^{\lambda_1}z_1,\dots,t^{\lambda_n}z_n)$ for all $t=
(t_1,\dots, t_d) \in T$.
Then,
\[
X_{A}= \overline{\mathcal{O}\bigl([1:\dots:1]\bigr)}=
\overline{\bigl\{[t^{\lambda_1}:\dots:t^{\lambda_n}]:\,t\in
(\Bbbk^{*})^d\bigr\}}
\subset \bP^{n-1}.
\]

By abuse of notation we also set $A\in \mathcal M_{d\times n}(\bZ)$ the
matrix with columns the weights $\lambda_i$.
In view of the  reductions made in Section~\ref{sec:prelim}
we assume without loss of generality that  the first row of $A$ is
$(1,\dots,1)$ and that the columns of $A$ span $\bZ^d$.

The homogeneous
ideal $I_A$ in $ \Bbbk[x_1,\dots,x_n]$ of the associated
projective toric variety $X_A$ is the binomial ideal (\cite{kn:stubook})
\[
I_A= \bigl\langle\, x^{a}-x^{b}:\;a, b \in
\mathbb{N}^{n}, \sum_{i=1}^{n}a_i\lambda_i = \sum_{i=1}^{n} b_i\lambda_i \,
\bigr\rangle .
\]
Thus, $X_A=\bigl\{[x] \in \bP^{n-1}:\;x^{a}=x^{b}, \;\forall\;a, b \in
\mathbb{N}^n\;\text{such that}\;Aa=Ab\bigr\},$ and it is easy to see that
\[
 X_{A}=\bigl\{[x] \in
\bP^{n-1}:\,x^{v^{+}}-x^{v^{-}}=0,\,\forall\,v\in
\mathcal R_A\bigr\},
\]
where $v_{i}^{+}=\max\{v_i, 0\}$, $v_{i}^{-}=-\min\{v_i, 0\}$ (and
so $v=v^{+}-v^{-}$).

For $p\in \bT^{n-1}$ we then have
\[
\begin{split}
m_{p}(X_A) = \overline{\mathcal{O}(p)} =
&  \bigl\{[x] \in
\bP^{n-1}\mathrel{:} p^{v^{-}}x^{v^{+}}-p^{v^{+}}x^{v^{-}}=0 ,\,\forall\,v \in
\mathcal R_A \bigr\}.
\end{split}
\]

\subsection{Characterization of self-duality in terms of the Gale dual
configuration} 

If $A$ is a non pyramidal configuration, then Theorem \ref{teoeqauto}
can be rephrased in terms of a geometric condition on the Gale dual of $A$.

\begin{definition} \label{def:Gale}
Let $A\in \mathcal M_{d\times n}(\bZ)$ with rank $d$. 
Let $\mathcal B_A=\{u_1,\dots, u_{n-d}\}\subset \bZ^n$ be a basis of $\mathcal R_A$.

We say that the matrix $B_A\in \mathcal M_{n\times (n-d)}(\bZ)$ with
columns the vectors $u_i$ is a {\em Gale dual matrix}\/ of $A$.
Let   $\mathcal G_A=\{b_1,\dots, b_{n}\}\subset \bZ^{n-d}$ be the configuration
of rows of $B_A$, that is $B_A=\left(\begin{smallmatrix}b_1\\
{\scriptstyle \vdots}\\ b_n \end{smallmatrix}\right)$  
(observe that we allow repeated elements).  The configuration $\mathcal G_A$
is called a {\em Gale dual}\/ configuration of $A$. Remark that
$\sum_{i=1}^n b_i =0$. 
\end{definition}

\begin{remark}
$(1)$ Since $\mathcal R_A$ is an affine invariant of the configuration
$A$, it follows that two affinely equivalent configuration share their
Gale dual configurations.

\noindent $(2)$ The configuration $A$ is non pyramidal if and only  $b_i\neq 0$
 for all $i=1,\dots, n$.

\noindent $(3)$ When $A$ is regular, $\mathcal R_A$ is the integer kernel
$\Ker_\bZ(A)$ of the matrix $A$.

\noindent $(4)$ For any Gale dual matrix of $A$, the morphism $\gamma:
\Bbbk^{n-d} 
\to \Bbbk^n$, $\gamma(s)=\bigl(\langle s, b_1
\rangle, \dots, \langle s, b_n \rangle\bigr)$ gives a parameterization
of $\mathcal R_{A,\bk}$, where we denote $\langle s, b_i \rangle = \sum_{j=1}^{n-d}
s_j b_{ij}$. 
\end{remark}

\begin{remark}\label{rmk:exproof}
By Theorem \ref{teoeqauto}, $X_A$ is self-dual if and only if
there exists $p_0 \in \bP
\bigl(\mathcal R_{A,\bk}\bigr) \cap \mathbb{T}^{n-1}$ such that $\bP
\bigl(\mathcal R_{A,\bk}\bigr) \subset \overline{\mathcal{O}(p_0)}$. 
By the remarks in subsection~\ref{rk:ector}, it follows that $X_A$ 
is self-dual if and only if for some such $p_0$
we have that $p_0^{v^-}w^{v^+} - p_0^{v^+}w^{v^-}=0$ for all $w\in \mathcal R_{A,\bk}$ and
$v \in \mathcal R_A$.

Assume $X_A$ is self-dual. Then, given any choice of Gale dual configuration,
 we deduce that
for all  $s \in \Bbbk^{n-d}\setminus\{0\}$ and  $j=1,\dots,n-d$, we have that
\[
p_{0}^{u_{j}^{-}}\bigl(\langle s, b_1
\rangle, \langle s, b_2 \rangle, \dots, \langle s, b_n \rangle
\bigr )^{u_{j}^{+}} \, = \, p_{0}^{u_{j}^{+}}\bigl(\langle s, b_1
\rangle, \langle s, b_2 \rangle, \dots, \langle s, b_n \rangle
\bigr )^{u_{j}^{-}}.
\]
for (some, or in fact all) $p_0 \in \bP
\bigl(\mathcal R_{A,\bk}\bigr) \cap \mathbb{T}^{n-1}$.

Since this gives an equality  in the polynomial ring $\Bbbk[s_{1}, \dots, s_{n-d}]$, 
both sides must have the same irreducible factors. But $\langle s,
b_i \rangle$ and $\langle s, b_k \rangle$ are associated
irreducible factors if and only if $b_i$ and $b_k$ are collinear
vectors. We deduce that for any line $L$ in $B$-space $\bZ^{n-d}$
and for all $j$,
\[
\sum \limits_{b_i \in L, b_{ij}
>0} b_{ij}= - \sum \limits_{b_i \in L, b_{ij}
<0} b_{ij}.
\]
Hence, $\sum \limits_{b_i \in L} b_{ij}=0$ for all
$j = 1, \dots, n-d$, or equivalently, $\sum \limits_{b_i \in L}b_i =0$.
\end{remark}

In fact, this last condition is not only necessary but also
sufficient. We give a proof of both implications using results about 
the tropicalization of the dual variety $X_A$ as described 
in \cite{kn:DFS05}.

First we recall that given a dual Gale configuration $\mathcal
G_A=\{b_1,\dots ,b_n\}$, and a subset $J\subset \{1, \dots, n\}$, 
the {\em flat}\/ $S_J$ of $\mathcal G_A$  associated to $J$ is the subset of
all the indices $i \in \{1, \dots, n\}$ such that $b_i$ belongs to the
subspace generated by $\{b_j\mathrel{:} j\in J\}$.

\begin{theorem}
 \label{teo:Bside}
Let $A\in \mathcal{M}_{d \times n}(\mathbb{Z})$ a non pyramidal 
configuration and $B_A$ a Gale dual for $A$ as in (\ref{def:Gale}). 
Then $X_{A}$ is self-dual if and only if for any line $L$ through the
origin in $\bZ^{n-d}$ we have that $\sum \limits_{b_{i} \in L} b_{i}=0$.
\end{theorem}

\begin{proof}
Since we are dealing with affine invariants, we can assume that $A$ is
a regular configuration. 
By Theorem~\ref{thm:autdim}, we know that $X_A$ is self-dual
if and only if $\dim X_A$ equals $\dim X_A^*$.  Given a vector
$v \in \bZ^n$, we define a new vector $\sigma(v) \in \{0,1\}^n$ by
$\sigma(v)_i = 0$ if $v_i \not =0$ and $\sigma(v)_i=1$ if $v_i=0$.

If follows from
\cite[Corollary~4.5]{kn:DFS05} that $\dim X_A = \dim X_A^*$  if and only if for any vector
$v \in \mathcal R_A$, the vector $(1, \dots, 1) - \sigma(v)$ lies in the row span $F$ of
the matrix $A$. But since we are assuming that $(1, \dots, 1) \in F$,
this is equivalent to the condition that $\sigma(v) \in F$.
 By duality, this is in turn equivalent to the fact that
for any $j=1, \dots, n-d$, the inner product 
\[
\langle \sigma(v), u_j \rangle \, = \,  \sum_{v_i=0} b_{ij} =0.
\]
That is to say, $X_A$ is self dual if and only if for any $v \in \mathcal R_A$ the sum $\sum_{v_i=0} b_i=0$.
But the sets $S$ of non zero coordinates of the vectors in the space of linear relations
$\mathcal R_A$ coincide with the flats of the Gale configuration $\mathcal G_A$. So, $X_A$
is self-dual if and only if for any flat $S \subset \{1, \dots, n\}$ the sum
$\sum_{i \in S} b_i=0$. It is clear that this happens if and only if the same condition
holds for all the one-dimensional flats, i.e. if for any line $L$ through the origin the sum
$\sum_{b_i \in L} b_i =0$.
\end{proof}

The assumption that  $A$ is a non 
pyramidal configuration in Theorem \ref{teo:Bside} is crucial, 
as  the following example shows.

\begin{example}
Let $A$ be a configuration such that $\mathcal R_A$ has rank $1$. Then
$\mathcal R_A$ is spanned by a single vector, whose coordinates add up
to $0$. So, the condition  
in Theorem~\ref{teo:Bside} that
the sum of the $b_i$ in this  line equals $0$ 
is satisfied. But by Corollary \ref{coro:noname} if $A$ is a pyramid,
then $X_A$ is not self-dual.
\end{example}

\subsection{Geometric characterization of self-dual configurations}
\label{sec:Aside}

In this paragraph we characterize the non pyramidal configurations
$A\subset \bZ^d$
whose Gale dual configurations are as in Theorem
\ref{teo:Bside}. We keep the assumptions  that  $\langle
A \rangle_\bZ= \bZ^d$ and that $A$ is non pyramidal.
We begin with some basic definitions about configurations.

\begin{definition}
Given $a = (a_1, \dots, a_n) \in \mathcal R_{A,\bk}$, we
call $\{i\mathrel{:} a_i\neq 0\}$ the {\em support}\/ of the relation and
denote $\supp(a)=\{i\mathrel{:}
a_i\neq 0\}$. We say that $\lambda_i$ belongs to the relation if $i\in
\supp(a)$.

Recall that any affine relation $a \in \mathcal R_{A, \bk}$ satisfies $\sum_i a_i=0$.
It is said that $a$ 
is a {\em circuit}\/  if there is no non trivial affine dependency relation with
support strictly contained in $\supp(a)$. In other words, a circuit is a minimal affine
dependency relation.
\end{definition}

\begin{remark}
\label{rem:circuitnumber}
Let $C$ be a circuit of a configuration $A$, and let $F$ be the
minimal face of $\operatorname{Conv(A)}$ containing $C$. If $d'$ denotes the dimension of the affine span
of $F$, then $C$ has at most $d'+2$ elements. 
\end{remark}

\begin{definition}
Two elements $b,b'$ of a configuration $B$ are {\em parallel}\/ if they generate
the same straight line through the origin. In particular, $b\neq 0$ and
$b'\neq 0$. The elements $b,b'$ are {\em antiparallel}\/ if they
are parallel and point into opposite directions.

Two elements $\lambda,\lambda'$ of a configuration $A$ are {\em
  coparallel}\/ if they 
belong exactly  to the same circuits.
\end{definition}

\begin{remark}
$(1)$ Coparallelism is an equivalence relation. We denote
by $\cc(\lambda)$ the coparallelism class of the element $\lambda\in
A$. 

\noindent (2) It is easy to see that $\lambda$ and $\lambda'$ are
coparallel if and only if they belong to the same affine dependency
relations. 

\noindent (3) The definition of coparallelism  can be extended to
pyramidal configurations as follows.  If $\lambda \in A$ is such that
it does not  belong to any dependency relation, then 
$\cc(\lambda)=\{\lambda\}$. Otherwise,  $\cc(\lambda)$ consists, as
above, of all elements of $A$ belonging to the same circuits as
$\lambda$. The condition that $A$ is not a pyramid is then equivalent to
the condition that $|\cc(\lambda)| \ge 2$ for all $\lambda \in A$.
 \end{remark}

\begin{lemma}
\label{lem:copara}
Let $\mathcal G_A = \{b_1, \dots, b_n\}$ be a Gale dual of $A$. Then $\lambda_i$
is coparallel to $\lambda_j$ if and only if
$b_{i}$ and $b_{j}$ are parallel elements of $\mathcal G_A$.
\end{lemma}

\begin{proof} Let $B_A$ the $(n \times (n-d))$-matrix with rows given by $\mathcal G_A$ as in
Definition~\ref{def:Gale}. As $A$ is not a pyramid, no row $b_i$ of $B_A$ is zero.
Any element $a \in \mathcal R_{A,\bk}$ is of the form $B_A \cdot m$,
for some $m \in \bk^{n-d}$. Then $\lambda_i$ is coparallel to $\lambda_j$ if and only if
for any nonzero $m \in \bk^{n-d}$ it holds that $\langle b_i, m \rangle \not=0$ precisely when
$\langle b_j, m \rangle \not=0$. It is clear that this happens if and
only if $b_i = \alpha b_j$ for 
a non zero constant $\alpha \in \bk$, that is, if and only if $b_i, b_j$ are parallel.
\end{proof}

\begin{definition}
Let $A=\{\lambda_1,\dots, \lambda_n\} \subset \mathbb Z^d$ be a
configuration. A subconfiguration $C' \subset A$ is called {\em facial} if there
exists a face $F$ of the convex hull $\operatorname{Conv}(A)\subset
\mathbb R^d$ of $A$ such that $C' = A \cap F$.

A subconfiguration   $C\subset A$ is a
{\em face complement}\/ if $A\setminus C$ is a facial subconfiguration
of $A$.
\end{definition}

\begin{remark}
\label{rem:deprel}
Let $A=\{\lambda_1,\dots, \lambda_n\} \subset \mathbb Z^d$ be a
configuration. A subconfiguration  $C=\{\lambda_{i_1},\dots,
\lambda_{i_h}\}\subset A$ is a
face complement if and only if there exists 
a dependency relation such that 
\[
\sum_{j=1}^h  r_{i_j} b_{i_j} =0 \quad , \quad r_{i_j} > 0.
\] 

Indeed,  a  dependency relation $\sum_{j=1}^h  r_{i_j} b_{i_j}
=0 $ with all 
$r_{i_j} > 0$  can be extended with zero coordinates 
to a relation $ r= (r_1, \dots, r_n)$ among all $b_i$'s. Thus, $r$
lies in the row 
space of $A$ and so there exists $\ell =(\ell_1, \dots, \ell_d)$ such 
that $r_i = \langle  
\ell, \lambda_i
\rangle$. It follows that  the linear form associated to $\ell$
vanishes on the complement of $C$, and 
all the points of $C$ lie in the same open half space delimited by the
kernel of $\ell$.
\end{remark}

\begin{lemma}
\label{lem:paracop}
Let $A=\{\lambda_1,\dots, \lambda_n\} \subset \mathbb Z^d$ be a
configuration. A coparallelism class  $C=\{\lambda_{i_1},\dots,
\lambda_{i_h}\}\subset A$ is a
face complement if and only 
if and only if there exist $j,k\in \{1,\dots ,h\}$ such that $b_{i_j}$ and
$b_{i_k}$ are antiparallel. 
\end{lemma}

\begin{proof}
If $C$ is a coparallelism class, we know by Lemma~\ref{lem:copara} that all
$b_{i_1}, \dots, b_{i_h}$ are parallel. It is then clear that a
dependency relation $r$ as in Remark \ref{rem:deprel} exists if 
and only if two of the vectors $b_{i_j}, b_{i_k}$ are antiparallel.
\end{proof}

\begin{definition}\label{def:pfc}
Let $A\subset \mathbb Z^d$ be a configuration and $C\subset A$ a face
complement. We say that $C$ is a {\em parallel face 
complement} if $C$ and  $A \setminus C$ lie in parallel hyperplanes. 

Note that in this case both $C$ and $A\setminus C$ are facial.
\end{definition}

\begin{example}
In Figure 1 below there are three configurations of $6$ lattice points in
$3$-dimen\-sional 
space (the $6$ vertices in each polytope). The $2$ vertices marked with big dots 
in each of the configurations define a coparallelism class $C$. In the first polytope
$(1)$, $C$ is not a face complement; in the second polytope $(2)$, $C$ is a face
complement but not a parallel face complement; in the third polytope $(3)$,
$C$ is a parallel face complement.  The characterization in
our next theorem proves that only the toric
variety corresponding to this last configuration is self-dual.
\begin{figure}[h]
\includegraphics{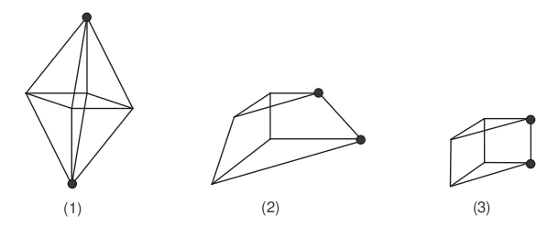}
\caption{Only configuration $(3)$ is self-dual.}
\end{figure}
\end{example}

It is straightforward to check that if $A_1, A_2$ are affinely
equivalent configurations and 
$\varphi$ is an affine linear map sending bijectively $A_1$ to $A_2$,
then $\varphi$ preserves 
coparallelism classes, faces and parallelism relations. Indeed, all
these notions can be read 
in a common Gale dual configuration. Moreover, we can translate
Theorem~\ref{teo:Bside} as follows.

\begin{theorem} \label{teo:Aside}
Let $A\subset \mathbb Z^d$ be a non pyramidal configuration.
The projective toric variety $X_A$ is self-dual if and only if any
coparallelism class of $A$ is a parallel face complement.
\end{theorem}

\begin{proof}
Let $\mathcal G_A$ be
a Gale dual of $A$ as in Definition~\ref{def:Gale}. By Lemma~\ref{lem:copara},
coparallelism classes $C =
\{ \lambda_{i_1} \dots, \lambda_{i_h}\}$ in $A$ are in correspondence with
parallel vectors $b_{i_1}, \dots, b_{i_h}$ in the dual space (i.e. lines  containing vectors of $\mathcal G_A$).
But now, $C$ is a parallel face complement if and only if there exists $\ell = (\ell_1, \dots, \ell_d)$
such that $\langle \ell, \lambda_i \rangle =0$ for all $\lambda_i\notin C$ and
$\langle \ell, \lambda_{i_j} \rangle =1$ for all $j=1, \dots, h$. Reciprocally, 
the sum of the vectors $\sum_{j=1}^h b_{i_j} =0$ implies the existence of such an $\ell$ as in
 Remark~\ref{rem:deprel}. The result now follows from
Theorem~\ref{teo:Bside}. 
\end{proof}

We have the following easy lemma.

\begin{lemma}\label{lem:ccpfc}
Assume that $A$ is a non pyramidal self-dual lattice configuration. Then, for any $\mu \in A$,
the coparallelism class $\cc(\mu)$ has at least two elements and it is a facial
subconfiguration of $A$.
\end{lemma}

\begin{proof}
It follows from Definition~\ref{def:pfc} that there exists a linear function $f$ taking
value $0$ on $A \setminus \cc(\mu)$ and value $1$ on $\cc(\mu)$. Then, $\cc(\mu)$ is
the facial subconfiguration of $A$ supported by the hyperplane $(f-1)=0$. If
$\cc{(\mu)} = \{\mu\}$, then by Theorem \ref{teo:Aside} $\{\mu\}$ is a
vertex, and hence $A$
would be a pyramid. It follows that  so $|\cc(\mu)|  \ge 2$, for any $\mu \in A$.
\end{proof}

We give in Lemma~\ref{lem:easy} $(2)$ an example of a self-dual lattice configuration $A$ which
contains an interior point of $\operatorname{Conv}(A)$. However, this cannot happen if
$X_A$ is not a hypersurface, as the following proposition shows.

\begin{proposition}
\label{prop:onlyvertex}
Let $A\subset \cX(T)$ be a configuration without repetitions such that
$X_A$ is self-dual, with $\operatorname{codim} X_A>1$. Then, the
interior of the 
convex hull $\operatorname{Conv}(A)$ does not contain elements of $A$ and for any
facial subconfiguration $C'$ of $A$, at most one point of $C'$ lies in the relative interior
of $\operatorname{Conv}(C')$.
\end{proposition}

\begin{proof}
Since $A=\{\lambda_1,\dots ,\lambda_n\}\subset \cX(T)$ has no repeated
elements,  it follows from Theorem \ref{teo:constaut+gde} that $A$ is
non pyramidal. Then, as $X_A$ is not a hypersurface, if follows from 
Remark \ref{rem:variado} that $n\geq d+3$, where $d$ is the dimension
of the affine span of $A$. 

Assume that there exists $\mu \in A$ belonging to 
the relative interior of an $s$-dimensional face $F$ of
$\operatorname{Conv}(A)$. Therefore, $\mu$ is a convex combination of the
vertices of $F$, and thus  $\cc(\mu)\subset F$. But  by Lemma~\ref{lem:ccpfc},
$\cc(\mu)$ is a facial subconfiguration of $A$, 
and thus a facial subconfiguration of $F \cap A$, which intersects the
relative interior of $F$. Then, $\cc(\mu) = F \cap A$. Let $\cc(\mu)=\{\mu,\lambda_1,\dots, \lambda_{r}\}$.
We claim that $\{\lambda_1,\dots, \lambda_{r}\}$ are affinely
independent and thus, $\cc(\mu)$ is a circuit. Indeed, for any
$i=1, \dots, r$, $\cc(\lambda_i) = \cc(\mu) = F \cap A$,  and so there cannot be any non-trivial
affine dependence relation involving only $\{\lambda_1,\dots, \lambda_{r}\}$.
In particular, $r=s+1$, $\{\lambda_1,\dots, \lambda_{s+1}\}$
are the vertices of $F$ and $\mu$ is the only  point  in $F \cap A$ belonging to
the relative interior of $\operatorname{Conv}(F)$.

Therefore, if the relative interior of $\operatorname{Conv}(A)$
contains one 
element $\mu\in A$, it follows that $A$ is a circuit, and hence
$n=d+2$, see Remark \ref{rem:circuitnumber}. That is, $X_A$ is a hypersurface. 
\end{proof}

\begin{example} \label{ex:intpointface}
Consider the self-dual configuration $A$ given by the columns of the matrix
\[
A=
\left(
\begin{matrix}
1 & 1 & 1& 0& 0& 0\\ 
0 & 0 & 0 & 1 & 1 & 1\\
0& 1 & 2& 0 & 0 & 0\\
0 & 0 & 0 & 0 & 1 & 2
\end{matrix}
\right).
\]
The associated toric variety has dimension $3$ in $\bP^5$, so it is not a hypersurface.
No point of $A = \operatorname{Conv}(A) \cap \bZ^4$ lies in the interior, but there are two facial
subconfigurations of $A$
(namely, the segments with vertices $\{(1,0,0,0), (1,0,2,0)\}$ and $\{(0,1,0,0), (0,1,0,2)\}$, respectively)
which do have a point of $A$ in their relative interior. Note that 
\[
X_A = \{(x_1, \dots, x_6) \in \bP^5 \, / \, x_2^2 - x_1 x_3 = x_4^2 - x_5 x_6\}
\]
 is not smooth. It is
a complete intersection but the four fixed points\\
 $(0,0,1,0,0,0), (1,0,0,0,0,0),(0,0,0,1,0,0), (0,0,0, 0,0,1)$ are not regular, as can be checked by
the drop in rank of the Jacobian matrix. This could be seen directly in the geometry of the configuration. 
The convex hull of $A$ is a simple polytope (in fact, it is a simplex) of dimension $3$ lying in the
hyperplane $H = \{(y_1, \dots, y_4) \in \bR^4\, / \, y_1 + y_2 =1\}$, but  fixing the origin at any
of the four vertices, the first lattice points in the
$3$ rays from that vertex do not form a basis of the lattice $H \cap \bZ^4$. Note that there is a 
splitting of the $4$-torus $T $ as a product of tori of dimension $2$ corresponding respectively
to the first three and last three weights in $A$.
\end{example}

We end this paragraph by showing another interesting combinatorial
property of configurations associated to self-dual toric varieties.

\begin{proposition}
\label{prop:passagetoface}
Let $A=\{\lambda_1,\dots, \lambda_n\}\subset \cX(T)$ be a  non
pyramidal configuration  such that $X_A$ is self 
dual and  let $D$ be an arbitrary non empty subset of $A$. Then, either $D$  is a  pyramidal configuration or
$X_D$ is self-dual and, moreover, $D$ is a facial subconfiguration of $A$.
\end{proposition}

\begin{proof}
Assume that $D=\{\lambda_1, \dots ,\lambda_s\}\subset A$ is non pyramidal, and consider $\mathcal
R_{D}\subset \mathbb Z^s$ . It is clear that 
$\mathcal R_{D}\times \{0\}\subset \mathcal R_{A}$. Hence, if
$\mathcal B_D$ is a basis of $\mathcal R_D$, then there exists a $\mathbb
Q$-basis of $\mathcal R_A\otimes \mathbb Q$ of the form $\mathcal
G_{D}\times \{0\} \cup \mathcal C$. Let $\mathcal B_A$ be a 
$\mathbb Z$-basis of $\mathcal R_A$, and  $\mathcal
G_A=\{b_1,\dots , b_n\}$ its associated Gale dual configuration. Then
there  exists an
invertible $\mathbb Q$-matrix $M$ such that
\[
B'=\left(\begin{array}{c|c} 
\ &\  \\
 B_D & C_1\\ 
\ &\ \\
\hline
\ &\ \\
0 & C_2\\
\ &\ 
\end{array}\right)=  \left(\begin{matrix}
b_1\\
\vdots\\
b_s\\
b_{s-1}\\
\vdots\\
b_n
\end{matrix}
\right)M
.
\]

Since $A$ is self-dual, it follows from Theorem \ref{teo:Bside} that the rows
$b_i$ are such that the sum of vectors $b_i$ in the same line through
the origin is zero. Hence, the matrix $B'$ satisfies the same
property. As $D$ is non pyramidal,  no row of $B_D$ is
zero. Therefore, $(B_D\,, C_1)$, and hence $B_D$,   also satisfy the
property that the sum of all its row vectors in a line through the
origin is equal to zero. Hence, 
$X_D$ is self-dual. Moreover, the sum of the  row vectors of $C_2$ is
zero, and it follows from Remark~\ref{rem:deprel} that $D$ is facial.
\end{proof}

\section{Families of self-dual projective varieties.} \label{sec:examples} 

In this section we use our previous results in order to obtain new
families of  projective
toric  varieties that are self-dual. In particular, we obtain many new
examples of non smooth self-dual projective varieties. We also identify 
all the smooth self-dual projective varieties of the form
$X_A$. We retrieve in this (toric) case Ein's result, without needing to rely on
Hartshorne's conjecture.

\subsection{Projective toric varieties associated to Lawrence
configurations}

\begin{definition}
\label{df:Law}
We say that  a configuration $A$ of $2n$ lattice points
is {\em Lawrence} if it is affinely equivalent to
a configuration whose associated matrix has the form
\begin{equation} 
\label{disp:Lawrence}
\left(
\begin{matrix} Id_n & Id_n\\ 0 & M \end{matrix} \right),
\end{equation}
where $Id_n$ denotes the $n \times n$ identity matrix.
Equivalently, $A$ is a Lawrence configuration if it is
affinely equivalent to a Cayley sum of $n$ subsets, each one containing
the vector $0$ and one of the column vectors of $M$.
\end{definition}

Lawrence configurations are a special case of \emph{Cayley
configurations} (see \cite{kn:cds01}). The Lawrence configuration 
associated to the matrix \eqref{disp:Lawrence} is the  Cayley
configuration of the two-point configurations consisting 
of the origin and one column vector of $M$. In the smooth case, Cayley
configuration of strictly equivalent polytopes correspond to toric fibrations
(see \cite{kn:DR03}).
 
It is straightforward to verify that if $A$ is Lawrence, then
\begin{enumerate}
\item[(i)] $\mathcal R_A=\left\{\left(\begin{smallmatrix} -v \\ \phantom{-} v
\end{smallmatrix}\right):\;v\in \Ker_\bZ(M)\right\}$.
\item[(ii)] $A$ is pyramidal  if and only if $M$ is pyramidal.
\end{enumerate}
We immediately deduce from Theorem~\ref{teo:Bside} the following result.

\begin{corollary}
\label{lauself-duales}
If $A$ is a non pyramidal Lawrence matrix then $X_A$ is self-dual.
\qed
\end{corollary}

\begin{example} \label{ex:segre}
The well known fact  that the Segre embedding of $\bP^1 \times \bP^{m-1}$ in
$\bP^{2m-1}$ is self-dual follows directly from
Corollary \ref{lauself-duales}, the image of
the Segre morphism
\[
\varphi (x,y) \, =
\,[y_0x_0:y_1x_0:\dots:y_mx_0:y_0x_1:y_1x_1:\dots:y_mx_1],
\]
where $x = [x_0, x_1], y = [y_0:y_1:\dots:y_m]$,
is a projective toric variety with associated matrix
\[
A=
\left(
\def\objectstyle{\scriptstyle}
\def\labelstyle{\scriptstyle}
\vcenter{
\xymatrix@=0pt
{
1\ar@{.}[rrr] &
 &  & 1 & 0\ar@{.}[rrr] &  &  & 0\\
0\ar@{.}[rrr] &  &  & 0 & 1\ar@{.}[rrr] &
 &  & 1 \\
 1\ar@{.}[dddrrr] & 0\ar@{.}[rr]\ar@{.}[ddrr] & & 0\ar@{.}[dd] & 1\ar@{.}[dddrrr]& 0\ar@{.}[rr]\ar@{.}[ddrr] &  & 0\ar@{.}[dd]\\ 
0\ar@{.}[dd]\ar@{.}[ddrr] &  &  &  & 0\ar@{.}[dd]\ar@{.}[ddrr] & &  & \\
 &  &  & 0 &  &  & & 0\\
0\ar@{.}[rr] &   & 0 & 1 & 0\ar@{.}[rr] &  & 0 & 1
}
}
\right)
\]

The sum of the first two rows equals the sum of the last $m$ rows. It
is easy to see that $A$ is 
affinely equivalent to the configuration $A'$ with associated matrix
\begin{equation} \label{eq:segre}
A'=
\left(
\begin{matrix}
Id_n &  Id_n\\
 \begin{smallmatrix}
0 &  \cdots  & 0
\end{smallmatrix} &
\begin{smallmatrix}
 1 &  \cdots  & 1
\end{smallmatrix}
\end{matrix}
\right).
\end{equation}
The matrix $A'$ is a non pyramidal Lawrence matrix, hence $X_{A'} =X_A$ is
self-dual.
\end{example}

We finish this paragraph by proving that Segre embeddings of
$\bP^1\times \bP^{m-1}$, $m\geq 2$ are the unique  smooth self-dual
projective toric varieties that are not a hypersurface.
We begin with an easy lemma which classifies all smooth hypersurfaces of the form $X_A$

\begin{lemma}\label{lem:easy}
Let $A$  be a lattice configuration such that $X_A$
is a smooth hypersurface. Then, $A$ is of one of the following forms:
\begin{itemize}
\item[(1)] $A$ consists of two equal points, and so $X_A = \{(1:1)\} = \{ (x_0:x_1) \in \bP^1 \, / \, x_0 - x_1=0\}$.  
\item[(2)] $A$ consists of three collinear points with one of them the mid point of the others, and so 
$X_A = \{(x_0:x_1:x_2) \in \bP^2\, / \, x_1^2 - x_0 x_2 =0\}$.
\item[(3)] $A$ consists of four points $a,b,c,d$ with $a+c=b+d$, and so $X_A =\{(x_0, x_1, x_2, x_3) \in \bP^3 \, /
\, x_0 x_3 - x_1 x_2 =0\}$ is the Segre embedding of $\bP^1 \times \bP^1$ in $\bP^3$.
\end{itemize}
\end{lemma}

\begin{proof}
When $X_A$ is a hypersurface, an equation for $X_A$ is given by
$b_A(x) = \prod_{b_i > 0} x_i^{b_i} - \prod_{b_i < 0} x_i^{-b_i}$,
where  the transpose of the row vector $(b_1, \dots, b_n)$ is a choice
of Gale dual of $A$. 
The cases $(1), (2)$ and $(3)$ in the statement correspond to the row vectors 
$(1,1), (1, -2,1)$ and $1,-1,-1,1)$, respectively (or any permutation of the coordinates), and it is 
straightforward to check that $X_A$ is smooth.  
It is easy to verify that in any other case, there exists a point
$x \in X_A$ where $b_A$ and all its partial derivatives vanish at $x$.
\end{proof}

We saw in Example~\ref{ex:intpointface} that a non-pyramidal self-dual
lattice configuration $A$ with ${\rm codim}(X_A) > 1$ can have a point
in the interior of a proper face. Moreover, more complicated
situations can happen: 

\begin{example}\label{ex:missingpts}
Consider the following dimension $3$ configuration 
$A \subset \bZ^4$, $A= \bigl\{(1,0,0,2), (1,0,0,0), (0,1,0,0), (0,1,0,2), 
(0,0,1,0),(0,0,1,1)\bigr\}$. Then, $\bZ A = \bZ^4$ and $X_A$ 
is self-dual because the following is a choice of Gale
dual $B \in \bZ^{6\times 2 }$:
\[
B=\left(\begin{array}{rr} 
1 & 0\\ 
-1 & 0\\
0 & 1 \\
0 & -1\\
2 & -2 \\
-2 & 2 
\end{array}\right).
\] 
All the points in $A$ are vertices of the polytope $P:=\operatorname{Conv(A)}$, but $A \not = P \cap \bZ^4$.
Indeed, there is a lattice point in the middle of each of the segments $[(1,0,0,2), (1,0,0,0)]$,
$[0,1,0,0), (0,1,0,2)]$, which are faces of $P$. It is clear that $X_A$ is not smooth (for instance
looking at the first lattice points in all the edges emanating from $(1,0,0,0)$), nor embedded by a
complete linear system.
\end{example}

However, the following result shows that when $X_A$ is smooth and self-dual, the situation is nicer.

\begin{lemma} \label{lem:descartando}
Let $A$ be a lattice configuration without repeated points such that $X_A$ is self-dual and \emph{smooth}.
Then, unless $X_A$ is  the quadratic rational normal curve in $(2)$ of Lemma~\ref{lem:easy}, no facial
subconfiguration $C \subseteq A$ contains a point of $A$ in the relative interior of $\operatorname{Conv}(C)$.
\end{lemma}

\begin{proof} 
Assume $A = \{ \lambda_1, \dots, \lambda_n\}$ 
has no repeated points and there exists $\mu \in A$ and a proper face $F$ of $\operatorname{Conv}(A)$
containing $\mu$ in its relative interior. Then, $F\cap A$ is not a
pyramid, and it follows from Proposition \ref{prop:passagetoface} that
$X_{F\cap A}$ is self-dual. Since  $X_{F \cap A}$ is also smooth,
Proposition~\ref{prop:onlyvertex} implies that  $X_{F\cap A}$ is a hypersurface. 
We deduce from Lemma~\ref{lem:easy} that $F \cap A$ has dimension one and consists  (up to reordering) of
 $3$ points  $\{ \lambda_1,
\lambda_2, \lambda_3\}$ with $ \lambda_1 + \lambda_3 = 2 \lambda_2$. 
We can choose a Gale dual $B$ of $A$ of the form:
\[
B=\left(\begin{array}{c|c} 
B_1 & C_1\\ 
\ &\ \\
\hline
\ &\ \\
0 & C_2\\
\ &\ 
\end{array}\right),
\]
with $B_1$ the $3 \times 1$ column vector with transpose $(1, -2,1)$. We see that the coparallelism class
of each $\lambda_i$ is contained in $F \cap A$ and no class can consist of a single element
because $A$ is not a pyramid. Therefore, $\cc(\lambda_i) = F \cap A, i = 1,2,3$; that is, any two of
the first $3$ rows of $B$ are linearly dependent. We can thus find another choice of Gale dual $B'$ of
$A$ of the form:
\[
B=\left(\begin{array}{c|c} 
B_1 & 0\\ 
\ &\ \\
\hline
\ &\ \\
0 & C_2\\
\ &\ 
\end{array}\right).
\]
Then, there is a splitting of the torus and $X_A$ cannot be smooth, with arguments similar to those
in Example~\ref{ex:intpointface}, because $A$ has no repeated points and so there is no linear
equation in the ideal $I_A$. 
\end{proof}

We now characterize the Segre embeddings 
$\bP^1 \times \bP^{m-1}$ in $\bP^{2m-1}$ from Example~\ref{ex:segre} in terms of the Gale dual configuration.

\begin{lemma}\label{lem:Bsegre}
A toric variety $X_A \subset \bP^{2m-1}$ is the Segre embedding of $\bP^1 \times \bP^{m-1}$ if and only if
any Gale dual $B \in \bZ^{2m \times r}$ of $A$ has the following form: $r = m-1$ and, up to reordering,
the rows of $b_1, \dots, b_{2m}$ of $B$ satisfy $\det( b_1, \dots, b_{m-1})=1, b_1+ \dots + b_m =0$ and
$b_{m+j} + b_j =0$, for all $j=1, \dots, m$.
\end{lemma}

\begin{proof}
It is clear that any Gale dual to the matrix $A'$ in~\eqref{eq:segre} is of this form. And it is also
straightforward to check that any matrix $B$ as in the statement is a Gale dual of this $A'$.
\end{proof}

We can now prove the complete characterization of smooth self-dual varieties $X_A$.

\begin{theorem}\label{th:regcase}
The only self-dual smooth non linear projective toric varieties equivariantly embedded are the toric hypersurfaces
described in $(2)$ and $(3)$ of Lemma~\ref{lem:easy} and the Segre embeddings
$\bP^1 \times \bP^{m-1}$ in $\bP^{2m-1}$ for $m \ge 3$. 
\end{theorem}

\begin{proof}
We proceed by induction in the codimension of $A$. By Lemma~\ref{lem:easy}, 
the result is true when $X_A$ is a hypersurface. Assume then that ${\rm codim}(X_A) > 1$. 
Now, by Lemma~\ref{lem:descartando}, we know that all the points in $A$ are vertices of 
$\operatorname{Conv}(A)$. Let $C$ be a coparallelism class and let
$D:= A\setminus C$. Then, $X_D$ is smooth and it is non
pyramidal. Indeed, we can choose a Gale dual $B$ of $A$ of the form:
\[
B=\left(\begin{array}{c|c} 
\begin{smallmatrix}b_{11}  \\ 
{\small{\vdots}}  \\
b_{r1}\end{smallmatrix}& 0 \\ 
\hline
\begin{smallmatrix}
b_{r+1,1}\\{\small{\vdots}}\\ b_{n1}\end{smallmatrix}& D_2\\
\end{array}\right),
\]
where $(b_{11},0),\dots, (b_{r1},0)$ correspond to the elements of $C$.  If
$D$ is a pyramid,  is is easy to show  that at least one row of $D_2$
must be zero, and it follows that the corresponding point of the
configuration belongs also to $C$, and thus is a contradiction. 

Hence, it follows from  Proposition~\ref{prop:passagetoface} 
that $X_D$  is self-dual with ${\rm codim}(X_D) = {\rm codim}(X_A) -1
< {\rm codim}(X_A)$ and no point of $D$ belongs to the relative interior of $\operatorname{Conv}(D)$. 
Therefore, by induction, $X_D$ is the Segre embedding of $\bP^1 
\times \bP^{m'-1}$ in $\bP^{2m'-1}$ for $m' \ge 2$ (including the hypersurface 
case $\bP^1 \times \bP^1$). In particular, $|D|= 2 m'$ is even.

Assume $C =\{ \mu_1, \dots, \mu_r\}$.
 Let $B_D \in \bZ^{2m'\times (m'-1)}$ be a choice of
Gale dual of $D$ as in Lemma~\ref{lem:Bsegre}, with rows $e'_1, \dots, 
e'_{m'}, - e'_1, \dots, - e'_{m'}$ with
$\{e'_1, \dots, e'_{m'-1}\}$ a basis of $\bZ^{m'-1}$ and $e'_1 + \dots + e'_m =0$. 
Add another integer affine relation
with coprime entries as the first column, to form a matrix $B'$ whose columns 
are a $\bQ$-basis of relations of $A$ of the form:
\[
B'=\left(\begin{array}{c|c} 
B_1 & 0\\ 
\ &\ \\
\hline
\ &\ \\
B_2& B_D\\
\ &\ 
\end{array}\right).
\]
Now, each coparallelism class of any $\mu \in D$ (with respect to $D$) 
has two elements when $m'>2$, and so it cannot
be ``broken'' when considering coparallelism classes in $A$, since it is not a pyramid. 
Then, via column operations we can assume that $B_2$ is of the form $B_2^t=(0, 
\dots, 0, a, 0, \dots, 0, -a),\,  (a \in \bZ_{\ge 0})$. In case $m'=2$,
then $B_D^t= (1, -1, -1,1)$ and the unique coparallelism class could be broken, but at most in two pieces with
two elements each, and again we have the same formulation for $B_2$.  In both cases, 
if $a=0$, then we
have a splitting, which implies that either there is a repeated point (if $B_1^t =( 1,-1)$) or $X_A$ is not smooth.
Then $a \ge 1$. Consider
the subconfiguration $E$ of $A$ obtained by forgetting 
the two columns corresponding to the rows $m'$ and $2m'$ of $B_D$. Since
the  vectors $b_i$ with complementary indices add up to zero, it follows that $E$ is facial and again, 
$X_E$ is smooth. We deduce that $a=1$ and $B_1^t = \pm (1, -1)$, which implies that $X_A$ is the Segre 
embedding of $\bP^1 \times \bP^{m'+1}$ in $\bP^{2m'+1}$.  
\end{proof}

\subsection{Non Lawrence families of examples}

We have the following obvious corollaries of Theorem~\ref{th:regcase}:

\begin{corollary}\label{cor:lardman}
Let $A \in \mathcal{M}_{d \times n}(\mathbb{Z})$ with maximal rank $d$
associated to a regular configuration of weights and let $X_A \subset \bP^{n-1}$ 
be the projective toric variety associated
to $A$. Assume $X_A$ is not a hypersurface, non linear, smooth and self-dual. Then, $n$ is
even. \qed
\end{corollary}

As the defect of the Segre embedding $X_m=\bP^1 \times \bP^{m-1}$ in $\bP^{2m-1}$ for any $m \ge 2$ equals
$2m-2 -m =m-2 = \dim X_m -2$, we recover for smooth varieties $X_A$ the following result, due to 
Landman (\cite{kn:Ein1}) for any projective smooth variety.

\begin{corollary}\label{cor:landman2}[Landman]
If $X_A \subset \bP^{n-1}$ is a non linear smooth
projective variety such that $\dim X < n-2$ with defect $k>0$, 
then $\dim X \equiv k\,(2)$.
\qed
\end{corollary}

We use  the previous corollaries together with Theorem~\ref{teo:Bside} 
to construct families of non regular self-dual varieties.

\begin{example}
\label{ej:familia}
Consider the families of matrices $\{A_{\alpha}\}, \{B_\alpha\}$ for $\alpha \in
\mathbb{Z}$, $\alpha \not=0,$ defined by:
\[ 
A_{\alpha}= \left(\begin{smallmatrix} 1 & 1 & 1 & 1 & \phantom{-} 1 & \phantom{-} 1 & 1 \\
1 & 1 & 1 & 1 & \phantom{-} 1 & \phantom{-} 0 & 0\\ 0 & 0 & 0 & 1 & \phantom{-} 1 & \phantom{-} 0 & 0
\\ 0 & 1 & 0 & \alpha & \phantom{-} 0 & -\alpha & 0\\ 0 & 0 & 1 & 0 & -\alpha & \phantom{-} 0 & \alpha
\end{smallmatrix} \right), \quad 
B_{\alpha}=\left(
\begin{smallmatrix} 
\phantom{\;} 2 \alpha & \phantom{-} 0 \\ -\alpha & \phantom{-} 0 \\ 
-\alpha  &\phantom{-} 0 \\ \phantom{-} 1 & \phantom{-} 1\\ -1 & -1 \\
\phantom{-} 0 & \phantom{-} 1 \\\phantom{-} 0 & -1 \end{smallmatrix}
\right).
\]

Clearly, $B_\alpha$ is   a choice of a Gale dual matrix of $A_\alpha$.

Observe that as $\alpha \neq 0$, the configuration $A_\alpha$ is not a pyramid 
and $\text{dim}(X_{A_{\alpha}})=4$. Moreover, it
is easy to show that if $\alpha\neq \alpha'$, then $X_{A_\alpha}$ and
$X_{A_{\alpha'}}$ are not isomorphic as embedded varieties because they have different
degrees. The degree of $X_{A_\alpha} $ is the normalized volume of the convex hull of the points
in the configuration $A_\alpha$ (\cite{kn:stubook}) and it can be computed easily in terms of the Gale
dual configuration. 

Since the conditions of Theorem \ref{teo:Bside} hold, it follows that
$X_{A_\alpha}$ is self-dual for all $\alpha\in  \mathbb Z, \alpha\not=0$. 
Moreover,  $n=7$ is odd and so we deduce from Corollary~\ref{cor:lardman}  that 
$X_{A_\alpha}$ is a singular variety.  The difference between its dimension and
its defect is $4-1=3 \not\equiv 0\; (2)$.
\end{example}

We can generalize Example \ref{ej:familia} in order to construct
families of non degenerate projective toric self-dual varieties of arbitrary
dimension greater than or equal to $3$ and of arbitrary codimension
greater than or equal to $2$.

\begin{example}
\emph{Families of self-dual varieties of any dimension $\geq 3$}.   Let any $r \geq 2$ and  $\alpha_1, \dots,
\alpha_r$ non zero integer numbers satisfying $\sum_{i=1}^r \alpha_i =0$.
Consider the planar lattice configuration 
\[
\mathcal G_\alpha = \bigl\{ (\alpha_1, 0), \dots, (\alpha_r,0),
(0,1), (0,-1), (1,1), (-1,-1)\bigr\}.
\]

Let $A$ be any lattice configuration with Gale dual $\mathcal G_\alpha$.
Then, $A$ is not a pyramid and the associated projective toric variety $X_A \subset \bP^{r+3}$  
is self-dual by Theorem~\ref{teo:Bside}, with dimension $\dim X_A = (r+ 4) - 2-1= r+1$.

When $r=2$, the dimension of $X_{A_\alpha}$ is
$3$. The case $\alpha_1, \alpha_2 = \pm 1$ corresponds to the Segre embedding
of $\bP^1 \times \bP^2$ in $\bP^5$. Already for $\alpha_1, \alpha_2 = \pm 2$,
the configuration $A_\alpha$ does not contain all the lattice points in its
convex hull. If we add those ``remaining'' points to the configuration, 
the associated toric variety is no longer self-dual.
\end{example}

\begin{example}
\emph{Families of self-dual varieties of any codimension $\geq 2$}.
Using the same ideas of the previous example, we can
construct pairs $(A,B)$ with $A$ a non
pyramidal configuration and $B$ its Gale dual satisfying the hypothesis of
Theorem~\ref{teo:Bside}, so that $X_A$ is self-dual, with arbitrary codimension  $m \geq 2$. 

For any $r \geq 2$ set $n = 2 m+ r$. As usual,  $e_1, \dots, e_m$ denotes the canonical basis in
$\bZ^m$. For any choice of non zero integers $\alpha_1, \dots, \alpha_r$
with $\sum_{i=1}^r \alpha_i =0$ consider the following lattice
configuration in $\bZ^m$: 
\[  
\mathcal G_\alpha:=  \bigl\{\alpha_1 e_1, \dots, \alpha_r e_1, e_2,
-e_2, \dots, e_m, -e_m,  
e_1 + \dots + e_m, - (e_1 + \dots + e_m)\bigr\}.
\]

For any lattice configuration $A_\alpha\subset \bZ^n$ with this Gale dual, $A_\alpha$ is
not a pyramid and its associated self-dual toric variety $X_{A_\alpha} \subset \bP^n$ has
dimension $m+r-1$ and codimension $m$.
\end{example}

\section{Strongly self-dual varieties} \label{sec:strong} 

We are interested now in characterizing a
particular interesting case of self-dual projective toric varieties.

\begin{definition} \label{def:strong}
Let $A$ be a  regular lattice configuration without repetitions. 
We say that the projective variety $X_A \subset \bP^{n-1}$ is {\em strongly self-dual}  if $X_A$
coincides with $X_A^{*}$  under the canonical identification between
$\bP^{n-1}$ and its dual projective space  as in Remark \ref{rk:identesp-espdual}.
\end{definition}

We deduce from Theorem~\ref{teoeqauto} the following characterization of strongly self-dual 
projective toric varieties of the form $X_A$.

\begin{proposition} \label{prop:strong}
Let $A$ be a  regular lattice configuration without repetitions. 
Then $X_A$ is strongly self-dual if and only if
$\bP\bigl(\mathcal R_{A,\bk}\bigr) \subset X_A$.
\end{proposition}

\begin{proof}
If $X_A$ is strongly self-dual, the containment 
$\bP\bigl(\mathcal R_{A,\bk}\bigr) \subset X_{A}^{*}$ implies
that the condition $\bP\bigl(\mathcal R_{A, \bk}\bigr) \subset X_A$ is necessary.

Assume that this condition holds and $A$ has no repetitions. As we already observed, 
Theorem~\ref{teo:constaut+gde}
implies that $A$ is not pyramidal. Then, it follows
from Theorem~\ref{teoeqauto} that for any $q \in \bP \bigl(\mathcal R_{A,\bk} \bigr) \cap \mathbb{T}^{n-1}
 \subset X_A \cap \mathbb{T}^{n-1}$, $m_{q}\bigl(X_A\bigr)= X_{A}^{*}$. 
But since $q \in \mathcal{O}\bigl([1:\dots:1]\bigr)$, we deduce that
$m_{q}\bigl(X_A\bigr)=\overline{\mathcal{O}(q)}=\overline{\mathcal{O}\bigl([1:\dots:1]\bigr)}=X_A$, 
that is $X_A^*=X_{A}$.
\end{proof}

Using the same notation of Theorem \ref{teo:Bside}, we have:

\begin{theorem}\label{teo:stself-dual}
Let $A$ be a non pyramidal regular lattice configuration $A$ of $n$ weights spanning $\bZ^d$
and let $B_A$ be a Gale dual of $A$. Then:

$X_A$ is strongly
self-dual $ \Leftrightarrow\,\,\left\{\begin{array}{l}\text{(a)}\,
\text{For any line}\, L \,\text{through the origin}\\ \text{we
have} \sum \limits_{b_{i} \in L} b_{i}=0.\\\text{(b)} \prod
\limits _{\begin{subarray}{c} j=1 \\ b_{ji} >0
\end{subarray}}^{n} b_{ji}^{b_{ji}}=\prod \limits_{\begin{subarray}{c} j=1 \\
b_{ji} < 0 \end{subarray}}^{n}
b_{ji}^{-b_{ji}},\; i=1,\dots,n-d.\end{array}\right.$
\end{theorem}

In the above statement, we use the convention that $0^0 =1$.

\begin{proof}
Assume that $X_A$ is strongly self-dual. Then (a) holds by Theorem~\ref{teo:Bside}. 
By Proposition~\ref{prop:strong}, we know that
$\bP \bigl( \mathcal R_{A,\bk}
\bigr) \cap \mathbb{T}^{n-1} \subset X_A \cap \mathbb{T}^{n-1}$,
and this last variety is cut out by the $(n-d)$ binomials
\[
\prod_{\begin{subarray}{c} j=1 \\ b_{ji} >0
\end{subarray}}^{n} x_j^{b_{ji}}=\prod_{\begin{subarray}{c} j=1 \\
b_{ji} < 0 \end{subarray}}^{n}
x_j^{-b_{ji}},\;\;\forall\;i=1,\dots,n-d.
\]

Then, we have the following equalities, for all $s \in \bk^{n-d}$:
\begin{equation} \label{eq:strong}
\prod_{\begin{subarray}{c} j=1 \\ b_{ji} >0
\end{subarray}}^{n} \langle s, b_j \rangle
^{b_{ji}}=\prod_{\begin{subarray}{c} j=1 \\
b_{ji} < 0 \end{subarray}}^{n} \langle s, b_j \rangle
^{-b_{ji}},\;\;\forall\;i=1,\dots,n-d.
\end{equation}
We get the conditions (b) by evaluating respectively at $s=e_1, \dots, e_{n-d}$.

Conversely, condition (a) implies  the equalities $(\ref{eq:strong})$ of the polynomials
in $s$ on both sides up to constant, as in
Remark~\ref{rmk:exproof}. Then, condition (b) ensures that this constant is $1$. Therefore,
$\bP \bigl( \mathcal R_{A,\bk}
\bigr) \cap \mathbb{T}^{n-1} \subset X_A \cap \mathbb{T}^{n-1}$, and so $X_A$ is strongly self-dual
 by Proposition~\ref{prop:strong}.
\end{proof}

\begin{example}
Consider the matrix $A=\left(\begin{smallmatrix}
1 & 0 & 0 & 0 & 0 & 0 & 0 & \phantom{-} 1 & \phantom{-} 1 \\
0 & 1 & 0 & 0 & 0 & 0 & 0 & \phantom{-} 1 & \phantom{-} 1\\
0 & 0 & 1 & 0 & 0 & 0 & 0 & \phantom{-} 2 & \phantom{-} 0\\
0 & 0 & 0 & 1 & 0 & 0 & 0 & \phantom{-} 0 & \phantom{-} 2\\
0 & 0 & 0 & 0 & 1 & 0 & 0 & -2 & -2\\
0 & 0 & 0 & 0 & 0 & 1 & 0 & -1 & \phantom{-} 0\\
0 & 0 & 0 & 0 & 0 & 0 & 1 & \phantom{-} 0 & -1\\
\end{smallmatrix} \right)$. Observe that $A$ is non pyramidal.
A Gale dual matrix $B_A$ for $A$ is given by the transpose of the matrix
 $\left(\begin{smallmatrix} -2 & -2 & -2&-2 & \phantom{-} 4 & \phantom{-} 1 & \phantom{-} 1 & \phantom{-} 
 1 & \phantom{-} 1 \\ \phantom{-} 1 & \phantom{-} 1 & \phantom{-} 2&
 \phantom{-} 0& -2 & -1 & \phantom{-} 0 &-1 & \phantom{-} 0  \end{smallmatrix}\right).$

Clearly,  $B_A$ satisfies the conditions of Theorem
\ref{teo:stself-dual} and hence $X_A$ is strongly self-dual.
But note that $A$ is not a Lawrence configuration.
\end{example}

We conclude this section with the complete characterization of strongly
self-dual varieties of type $X_A$, with  $A$ a  non pyramidal Lawrence matrix.

\begin{theorem}\label{teo:caracLawself-dual}
Let $A$ be a non
pyramidal Lawrence configuration consisting of $2n$ points in
$\bZ^{n+d}$, as in $(\ref{disp:Lawrence})$. Then $X_A$ is strongly self-dual if
and only if there 
exists a subset $I$ of rows of the lower matrix $M= (m_{jk}) $ such that $\sum_{j\in I} m_{jk}$ is an odd
number for all $k=1, \dots, n$. 
\end{theorem}

\begin{proof}
By Corollary \ref{lauself-duales},  $X_A$ is self-dual for any non pyramidal Lawrence 
configuration $A$. Thus, $X_A$ is strongly self-dual if and only if conditions (b) in
Theorem~\ref{teo:stself-dual} are satisfied. 
If $\mathcal G_M =\{c_1, \dots, c_n\} \subset \bZ^{n-d}$ is a Gale dual configuration for
$M$, then $\{ -c_1, \dots, -c_n, c_1, \dots, c_n\}$ defines a Gale dual
configuration for $A$. Conditions (b) are then equivalent in this case to
the equalities
\[
(-1)^{\sum_{j=1}^n  c_{ji}} \, = \, 0, \quad i =1, \dots, n-d.
\]
This is in turn equivalent to the condition that for all $v \in \mathcal{R}_{M}$,
the sum $\sum_{j=1}^n v_j \equiv 0 \; (2)$. But this is equivalent to the fact that
the vector $(1, \dots, 1)$ lies in the row span of $M$ when we reduce all its entries
modulo $2$. Denoting classes in $\bZ_2$ with an over-line, this condition means
that there exist $\alpha_1, \dots, \alpha_d \in \bZ_2 = \{0,1\}$ such that
\[
(1, \dots, 1) = \sum_{i=1}^{d} \alpha_{i}(\overline{m_{i1}}, \dots,
\overline{m_{in}})
= \sum_{\alpha_i =1} (\overline{m_{i1}}, \dots,
\overline{m_{in}}).
\]

It suffices to call $I = \bigl\{ i \in \{1, \dots, d\} \, : \, \alpha_i =1\bigr\}$.
\end{proof}

\begin{example}
The Segre embeddings in Example~\ref{ex:segre}
have associated Law\-rence matrices as in (\ref{eq:segre}), where $M$
is a matrix with a single row of with all entries equal to $1$. They clearly
satisfy the hypotheses of Theorem~\ref{teo:caracLawself-dual}. Then,
for any $m >1$, the Segre embedding of $\bP^1 \times \bP^{m-1}$ is
a strongly-self dual projective toric variety.
\end{example}

\bigskip

\bigskip

\bigskip

\begin{small}

\begin{tabular}{l}
M.~Bourel \\
Instituto de Matem\'atica y Estad\'{\i}stica\\
Facultad de Ingenier\'{\i}a\\
 Universidad de la Rep\'ublica\\
Julio Herrera y Reissig 565\\
11300 Montevideo, Uruguay\\
\email{mbourel@fing.edu.uy}
\end{tabular} 
\begin{tabular}{l}
A.~Dickenstein\\
Departamento de Matem\'atica\\
FCEN, Universidad de Buenos Aires \\
Ciudad Universitaria, Pab.I \\
(1428) Buenos Aires, Argentina\\
\email{alidick@dm.uba.ar}
\end{tabular} 

\bigskip

\medskip

\begin{tabular}{l}
A.~Rittatore\\
Facultad de Ciencias\\
Universidad de la Rep\'ublica\\
Igu\'a 4225\\
11400 Montevideo, Uruguay\\
\email{alvaro@cmat.edu.uy}
\end{tabular}
\end{small}

\begin{thebibliography}{100}

\bibitem{kn:Bor} A.~Borel, \emph{Linear Algebraic Groups}, second
edition. Graduate Texts in Mathematics, \textbf{126}. New York,
Springer, 1991.

\bibitem{kn:bps} D.~Bayer, S.~Popescu and B.~Sturmfels,
\emph{Syzygies of Unimodular Lawrence Ideals}, J. Reine Angew. Math.
534  (2001), 169--186.


\bibitem{kn:DRC05} C.~Casagrande, S.~Di Rocco,  \emph{Projective Q-factorial
toric varieties covered by lines}, Commun. Contemp. Math.  10  (2008),
no. 3, 363--389. arXiv:math.AG/0512385.

\bibitem{kn:cds01} E.~Cattani, A.~Dickenstein,  B.~Sturmfels, \emph{Rational
Hypergeometric Functions}, Compositio Math.  128  (2001),  no. 2, 217--239.

\bibitem{kn:cls} D.~Cox, J.~Little, H.~Schenk, \emph{Toric Varieties},
  Graduate Studies in  Mathematics, Amer. Math. Soc., Providence, RI,
  to appear.
 
\bibitem{kn:CC05} R.~Curran, E.~Cattani,  \emph{Restriction of A-
    Discriminants and 
dual defect toric varieties}, J. Symbolic Comput. 42  (2007),
no. 1-2, 115--135.

\bibitem{kn:DS02} A.~Dickenstein, B.~Sturmfels, \emph{Elimination theory in
codimension two}, Journal of Symbolic Computation, 34 (2002),
119-135.

\bibitem{kn:DFS05} A.~Dickenstein, E.M.~Feichtner, B.~Sturmfels,
\emph{Tropical Discriminants},  J. Amer. Math. Soc.  20  (2007),
no. 4, 1111--1133. 


\bibitem{kn:DR03} S.~Di Rocco, \emph{Projective duality of toric
manifolds and defect polytopes}.  Proc. London Math. Soc. (3)  93
(2006),  no. 1, 85--104. 



\bibitem{kn:Ein1} L.~Ein, \emph{Varieties with small dual varieties}, I,
Invent. Math. 86 (1986), no. 1, 63--74.

\bibitem{kn:Ein2} L.~Ein, \emph{Varieties with small dual varieties}, II, Duke.
Math. J. 52 (1985), no. 4, 895--907.

\bibitem{kn:Gel} I.M.~Gelfand; M.M.~Kapranov; A.V.~Zelevinsky,
\emph{Discriminants, resultants, and multidimensional
determinants}. Boston, Birkh\"{a}user, 1994.

\bibitem{kn:Gel2}  I.M.~Gelfand; M.M.~Kapranov y A.V.~Zelevinsky,
\emph{Hypergeometric functions 
and toral manifolds}. Funct. Anal. Appl. 23 (1989), no. 2, 94--106.


\bibitem{kn:Pop02} V.L.~Popov, \emph{Self-dual algebraic varieties and nilpotents
orbits}. Proceedings of the International Colloquium on Algebra,
Arithmetic and Geometry, Mumbai, 2000, Tata Inst. Fund. Research,
Narosa Publ. House, 2002, 509--533.

\bibitem{kn:PT04} V.L.~Popov; E.A.~Tevelev, \emph{Self-dual projective algebraic
varieties associated with symmetric spaces}. In: Popov,
Vladimir L. (ed.), Algebraic  groups and algebraic
varieties, Encyclopaedia of Mathematical Sciences, Vol. 132, Invariant Theory and
Algebraic Transformation Groups Vol. III, Springer Verlag, 131--167 (2004). 


\bibitem{kn:stubook} B.~Sturmfels, \emph{Gr\"obner bases and convex polytopes},
University Lecture Series, 8. American Mathematical Society,
Providence, RI, 1996. xii+162 pp.. 

\bibitem{kn:tev05} E.~Tevelev, \emph{Projective duality and homogeneous
spaces}, Encyclopaedia of Mathematical Sciences, 133. Invariant
Theory and Algebraic Transformation Groups 4. New York, Springer,
2005.

\bibitem{kn:zak}  F.~Zak, \emph{Projection of algebraic varieties},
  (Russian)  Mat. Sb. (N.S.)  116(158)  (1981), no. 4, 593--602,
  608. English translation: Math. Sbornik {\bf 44} (1983), 535--544.  

\bibitem{kn:Zi} G.~Ziegler, \emph{Lectures on polytopes}, GTM,
  152. Springer-Verlag, New York, 1995. x+370 pp.
\end{thebibliography}
\end{document}